\input{style/arxiv-ba.cfg}
\documentclass[ba,linksfromyear,dvips,preprint]{imsart}
\usepackage{natbib}
\usepackage{amsthm,amsmath}
\usepackage{amssymb}
\usepackage{vtexurl}
\usepackage{amssymb}
\usepackage{amsthm}
\usepackage{graphicx}

\pubyear{2015}
\volume{10}
\issue{1}
\firstpage{139}
\lastpage{170}
\doi{10.1214/14-BA890}% Updated by VTEXPTS2LaTeX.exe, 22.01.2015 11:09

\startlocaldefs

\newcommand{\Be}{\mbox{{\bf Be}}}

\newcommand{\B}{\mbox{{\bf B}}}

\newcommand{\corr}{\mathbf{Corr}}
\newcommand{\nd}{\noindent}
\newcommand{\logit}{\mbox{logit}}

\newcommand{\X}{\mathcal{X}}
\newcommand{\x}{\mathbf{x}}

\theoremstyle{plain}
\newtheorem{thm}{Theorem}[section]

\newtheorem{prop}[thm]{Proposition}

\theoremstyle{definition}

\theoremstyle{remark}

\endlocaldefs

\begin{document}

\begin{frontmatter}
\title{Objective Bayesian Inference for Bilateral Data}
\runtitle{Objective Bayesian Inference for Bilateral Data}

\begin{aug}
\author{\fnms{Cyr Emile} \snm{M'lan}\thanksref{t1}}%,
\and
\author{\fnms{Ming-Hui} \snm{Chen}\thanksref{t2}}

\runauthor{Cyr Emile M'lan and Ming-Hui Chen}

% \affiliation{}

\thankstext{t1}{William E. Wecker Associates, Inc., 270 E Simpson
Avenue, Jackson, WY 83001, USA,\hfill\break {cmlan@hotmail.com}}
\thankstext{t2}{Department of Statistics, University of Connecticut,
Storrs, CT 06269, USA,\hfill\break {ming-hui.chen@uconn.edu}}

\end{aug}

%% Abstract %%
%
\begin{abstract}
This paper presents three objective Bayesian methods for analyzing
bilateral data under Dallal's
model and the saturated model. Three parameters are of interest,
namely, the risk difference,
the risk ratio, and the odds ratio. We derive Jeffreys' prior and
Bernardo's reference prior associated
with the three parameters that characterize Dallal's model. We derive
the functional forms of the
posterior distributions of the risk difference and the risk ratio and
discuss how to sample
from their posterior distributions. We demonstrate the use of the
proposed methodology with two real data examples.
We also investigate small, moderate, and large sample properties of the
proposed methodology and the frequentist
counterpart via simulations.
\end{abstract}

%% Keywords %%
%

% \begin{keyword}[class=MSC]
% \kwd[Primary ]{}
% \kwd[; secondary ]{}
% \end{keyword}

\end{frontmatter}

%% Mainmatter %%

\section{Introduction}\label{sec1}

Bilateral data arise in medicine when a group of randomly chosen
patients with a condition receive a
new treatment, for example, surgery, on paired body parts within the
same individual (eyes, ears, breasts, arms,
hands, knees, legs, or feet), while another group of patients with this
condition receive a control treatment, for example, the
currently accepted medical treatment. The investigator records paired
Bernoulli outcomes about a particular
characteristic, for example, absence of the condition, that are then
grouped into one of three categories.
(i) The two body parts are cured, recorded as $(1,1)$. The counts of
patients with this characteristic from the control and
treatment groups are denoted by $m_{20}$ and $m_{21}$. (ii) One of the
two body parts is cured while the other remains
diseased, recorded as $(1,0)$ or $(0,1)$. However, these $(1,0)$ and
$(0,1)$ outcomes are thrown away and this detailed
information is no longer available. Only the counts of patients from
the control and treatment groups in that category,
$m_{10}$ and $m_{11}$ are available. (iii) Neither of the two body
parts are cured, denoted
by $(0,0)$. The counts of patients in this category from the control
and treatment groups are denoted by $m_{00}$ and $m_{01}$.
The data, denoted by $D$, can be summarized into a $3\times2$
contingency table (see Table \ref{tab1sec1})
where the trinomial counts $(m_{01}, m_{11}, m_{21})$ for the treatment
group and $(m_{00}, m_{10}, m_{20})$ for the
control group are the cell entries. Such data are very common in ophthalmologic,
orthopaedic and otolaryngologic studies. Twin studies are also a
familiar source
of bilateral data. The goal of such clinical trials is to quantify the
benefit of treatment over placebo.

\begin{table}[h]\caption{Data structure for bilateral data along with
their corresponding trinomial
probabilities. $p_{00}+p_{10}+p_{20}=1=p_{01}+p_{11}+p_{21}$,
$m_{00}+m_{10}+m_{20}=m_{+0}$
and $m_{01}+m_{11}+m_{21}=m_{+1}$ are fixed by design.}\label{tab1sec1}
\bigskip\centering
\begin{tabular}{lllll}\hline
& \multicolumn{4}{c}{Group} \\\cline{2-5}
Numbers of & & & & \\
cured organs & Treatment & & & Control\\\hline
0 & $m_{01} \; (p_{01})$ & & & $m_{00} \; (p_{00})$\\
1 & $m_{11} \; (p_{11})$ & & & $m_{10} \; (p_{10})$\\
2 & $m_{21} \; (p_{21})$ & & & $m_{20} \;(p_{20})$\\
Total & $m_{+1}$ &&& $m_{+0}$\\\hline
\end{tabular}
\end{table}
The main parameter in bilateral models is the difference between the
proportion of body parts
$(\lambda_0,\lambda_1)$ with the characteristic of interest in the two
groups, $\Delta=\lambda_1-\lambda_0$.
However, the dependency between paired observations cannot be ignored.
\citet{Rosner1982}, \citet{Morris1993}, and \citet{TangTangRosner2006}
%Tang, and Rosner (2006)
discuss the consequences of ignoring this dependency.

Five models have been proposed for bilateral data: Rosner's model,
Dallal's model, the equal
correlation model, the independence model, and the full or saturated
model. Of these five
models, the most extensively studied is Rosner's model. \citet
{TangTangQiu2008} %Tang, Tang and Qiu (2008)
present several test statistics for the equality of $\lambda_0$ and
$\lambda_1$ under Rosner's model.
\citet{QiuTangTang2009} % Qiu, Tang and Tang (2009)
consider the problem of sample size calculations under Rosner's model.
\citet{TangQiuTangPei2011} % Tang et al. (2011)
discuss various techniques to construct asymptotic confidence intervals
for~$\Delta$ under
Rosner's model and evaluate the performance of these via empirical studies.
\citet{TangPeiWongLi2010} %Tang et al. (2010)
is perhaps the only paper that discusses all five models simultaneously
to determine which one
provides a better fit to the data.
\citet{PeiTangWongGuo2010} % Pei et al. (2010)
present asymptotic confidence intervals
under the equal correlation model and evaluate their performance.
\citet{PeiTangGuo2008} % Pei, Tang and Guo (2008)
present a test for the equality of $\lambda_0$ and $\lambda_1$ when
unilateral and bilateral data are
combined under the equal correlation model.

Among these five models, the independence model is rarely used in
practice in the context of bilateral data.
The objective priors for the saturated model have already been
investigated in the literature.
The challenge with Rosner's model is that it is difficult to justify
this model from a biological point of view.
Only the equal correlation model and the equal conditional probability
model (Dallal's model)
have sound statistical foundations and biological interpretation.
Indeed, one can always characterize a bivariate discrete distribution
for two binary random variables by making assumptions
about (i) their marginal distributions and the correlation they share
or (ii) their marginal distributions and the two conditional
distributions they share.
For Rosner's model and the equal correlation model, Jeffreys' prior and
Bernardo's reference prior
are too complex to be of much practical use. Dallal's model is the only
model for which we could derive useful closed-form expressions
for various objective priors. For these reasons, we solely focus on
Dallal's model in this
paper.\vadjust{\eject}

Bayesian methods have gained incredible popularity in recent years both
in the theory and practice
of statistics. Under non-informative priors, Bayesian inferences yield
results similar to that obtained
under the frequentist paradigm. % This is the case with the example
%discussed in Section \ref{sub2sec6}.
As we will show in an example presented in Section \ref{sub2sec6},
Bayesian inference for Dallal's model yields similar findings to the
corresponding frequentist analysis.
Bayesian methods, however, do not rely on the normal approximation to
carry out statistical inferences, which is an advantage
over the frequentist methods.
It is not unusual to encounter $3\times2$ bilateral data where one or
more cells have sparse data,
thereby preventing the use of the usual normal approximation that
underlies frequentist inferences
developed in Appendix 1 of the Supplementary Web Materials. Also, in
some sparse $3\times2$
bilateral data, frequentist estimates lie on the boundary of the
parameter space, which are not permitted
by design. In such situations, one cannot compute the confidence
intervals for some parameters or
carry out tests of hypotheses about some of the parameters. We provide
one such example in Section
\ref{sub1sec6}. Bayesian methods provide a simpler way to analyze such
sparse bilateral data.
Another benefit of Bayesian inference for $3\times2$ bilateral data is
its ability to handle
the nuisance parameter in Dallal's model which complicates frequentist analyses.

There are no existing Bayesian methods for bilateral data in the
literature yet. We present
objective Bayesian inferences for three parameters of interest: the
risk difference, the risk ratio, and
the odds ratio. In Section 2, we present Dallal's reduced model along
with Dallal's full model.
Section 3 is dedicated to the derivation of the objective Bayesian
modeling of bilateral data.
We focus primarily on deriving Jeffreys' prior and Bernardo's reference
prior. We then discuss the Bayes
factor in the context of hypothesis testing as well as a simulation
scheme for the joint posterior distribution.
Section 4 presents results of an empirical comparison between Bayesian
methods and frequentist methods.
Section 5 presents two illustrative case studies. In Section 6, we
present two families of Bayesian prior distributions,
including Jeffreys' or Bernardo's reference priors as special cases.
Section 7 concludes this paper.

%%----------------------------------------------------------------------------------------------
% Section 2
%%----------------------------------------------------------------------------------------------

\section{Dallal's Dependence Model}\label{sec2}

Dallal's model was presented for the first time as model 2 in \citet
{Dallal1988}. %Dallal (1988).
Let ``$i=1$" and ``$i=0$" denote the treatment group and the control
group, respectively.
Denote by $Z_{ijk}$
a binary variable such that $Z_{ijk}=1$ if the $k$th site ($k$th body
part) of the $j$th subject in the $i$th
treatment group is free of disease at the end of the study and 0
otherwise for $i=0,1$, $j=1,\dots, m_{+i}$, and $k=1,2$.
Dallal's (reduced) model is characterized by the following assumptions:
\begin{description}
\item[Assumption 1:] $P\big(Z_{ijk}=1\big) = \lambda_i$ for $i=0,1$
with $0<\lambda_i<1$.
\item[Assumption 2:] $P\big(Z_{ijk}=1 \mid Z_{ij(3-k)}=1\big) = 1-\gamma
$ with $0<\gamma<1$.
\end{description}

Assumption 2 states that the conditional probability of an occurrence
of a particular
characteristic at one site given an occurrence of that characteristic
at the other\vadjust{\eject} site to be the same in the two treatment
groups. This statement is relaxed and replaced by
$P\big(Z_{ijk}=1 \mid Z_{ij(3-k)}=1\big) = 1-\gamma_i$ with $0<\gamma
_i<1,\,i=0,1$
in the full or saturated model. That is, two conditional probability
statements are made, one for the
treatment group and the other for the control group. We also refer to
this saturated model as Dallal's saturated model.
However, the full model has one more parameter than the reduced model.

Let $m_{hi}$ be the number of subjects in the $i$th group with exactly
$h$ site(s) cured and
$p_{hi}$ be the success probability associated with $m_{hi}$ for
$h=0,1,2$ and $i=0,1$. The two group
total sample sizes are denoted by $m_{+1}$ for the treatment group and
$m_{+0}$ for the control group
and these are assumed fixed by design. Hence, $(m_{0i}, m_{1i},
m_{2i})$ follows the trinomial distribution
with total number of trials $m_{+i}$ and probability parameter vector
$(p_{0i}, p_{1i}, p_{2i})$ for $i=0,1$ such as
\[
p_{0i} = 1 - (1+\gamma)\lambda_i, \;\;
p_{1i} = 2\gamma\lambda_i, \;\; \mbox{and} \;\;
p_{2i} = (1-\gamma)\lambda_i.
\]
Dallal's model also implies that the correlation coefficients between
the $Z_{ijk}$ variables take the form:
\begin{eqnarray}
\rho_i = \corr\big(Z_{ijk}, Z_{ij(3-k)}\big) = 1- \dfrac{\gamma
}{1-\lambda_i},\qquad i=0,1.\label{eq1sec2}
\end{eqnarray}
In fact, the conditional probability assumption in Dallal's model can
be replaced by the statement in
\eqref{eq1sec2} about the correlation coefficient. The correlation
coefficient, $\rho_i$, takes both positive
and negative values over the entire range $(-1,1)$. The excess risk is
defined as
\[
\delta_i = P(Z_{ijk}=1 \mid Z_{ij(3-k)}=1) - P(Z_{ijk}=1 ) = 1-\gamma-
\lambda_i
\]
for $i=0,1$.

The main parameter of interest in this investigation is the risk
difference, $\Delta=\lambda_1 - \lambda_0$,
and, therefore, $\gamma$ can be viewed as a nuisance parameter. The
risk ratio, $R =\dfrac{\lambda_1}{\lambda_0}$
and the odds ratio, $\psi= \dfrac{\lambda_1(1-\lambda_0)}{(1-\lambda
_1)\lambda_0}$ can also be of interest.
Another parameter of interest is the difference of excess risks in both
the treatment and the control groups,
$\delta_0-\delta_1$, which is equal to $\Delta$ under Dallal's reduced model.
To date, the risk ratio, the odds ratio and the difference of excess
risks have never been discussed in the bilateral
data literature. The first two parameters add another dimension to the
utility of bilateral data so that they
can be collected under either a prospective study paradigm or a
retrospective observational study paradigm,
allowing for more applications than those under the current clinical setting.
Although the addition of these new
parameters poses more challenges for carrying out frequentist
inference, no additional
work is required in the Bayesian framework.\vadjust{\eject}

What makes the inferential process challenging in Dallal's model is
that one deals with a
constrained parameter space. Indeed, the parameter space is
\begin{eqnarray*}
\Omega= \bigg\{ \big(\gamma, \lambda_0, \lambda_1\big): \; & & 0
< \gamma< 1 \mbox{ if } 0 < \max(\lambda_0, \lambda_1) \leq\frac
{1}{2};\nonumber\\
\; & & 0 <\gamma< \dfrac{1}{\max(\lambda_0, \lambda_1)}-1 \mbox{
if } \max(\lambda_0, \lambda_1)>\frac{1}{2} \bigg\}\;\;
\end{eqnarray*}
or equivalently
\begin{eqnarray}\label{eq1sec3}
\Omega= \bigg\{ \big(\gamma, \lambda_0, \lambda_1 \big): \;
& & 0<\gamma<1\;\mbox{and}\; 0 < \lambda_0, \lambda_1 <
\dfrac{1}{1+\gamma} \bigg\}.
\end{eqnarray}
We adopt the second representation in the sequel.
%Note that the parameter values $\lambda_0=\lambda_1=\gamma=
%the special case of independence.

The likelihood function can be expressed as
\begin{eqnarray*}
L(\gamma, \lambda_0, \lambda_1) & = & \prod_{i=0}^1 {m_{+i} \choose
{m_{0i}, m_{1i}, m_{2i}}} \big[1-(1+\gamma)\lambda_i\big]^{m_{0i}}
(2\gamma\lambda_i)^{m_{1i}} \big[(1-\gamma)\lambda_i\big
]^{m_{2i}}\nonumber\\
& \propto& \gamma^{m_{10}+ m_{11}}
(1-\gamma)^{m_{20}+m_{21}}\lambda_0^{m_{10} + m_{20}} \lambda_1^{m_{11}
+ m_{21}}\;
\big[1-(1+\gamma)\lambda_0\big]^{m_{00}}\nonumber\\
& & \times\big[1-(1+\gamma)\lambda_1\big]^{m_{01}},\qquad\qquad\qquad
(\gamma,\lambda_0,\lambda_1)\in\Omega.
\end{eqnarray*}
In Appendix 1 of the Supplementary Web Appendix, we develop for the
first time large-sample frequentist
inferences for the risk difference, the risk ratio and the odds ratio
under Dallal's model.

%%----------------------------------------------------------------------------------------------
% Section 3
%%----------------------------------------------------------------------------------------------

\section{Bayesian Analysis}\label{sec4}

A key component in Bayesian analysis is the choice of the prior
distribution. Traditionally,
Bayesians have turned to conjugate priors. However, the concept of
conjugate priors is
only universal in standard univariate problems. In addition, when
little or no prior information
is available, conjugate priors become subjective. Thus, non-informative
or objective priors
are more widely accepted. The uniform distribution over the parameter
space is
an obvious non-informative prior. Jeffreys' prior and Bernardo's
reference prior are
alternative choices that are invariant under any parameterization or a
larger class of
parameterizations. Another alternative choice is to use a joint prior
that is the
compromise between an informative prior and a non-informative prior
\citep{SunBerger1998}. %(Berger and Sun, 1988).

We discuss four types of priors: the uniform prior, Jeffreys' prior,
Bernardo's prior,
and Sun and Berger's reference prior in light of partial information.
We derive the posterior distributions of $\Delta$ and $R$ as well as
$P\big(\Delta>\Delta_0\,|\,D\big)$ and $P\big(R>R_0\,|\,D\big)$ in
Appendix 4. %\ref{sec4app}.
Although we do not provide the posterior distribution of $\psi$ (more
complex), we discuss in Section
\ref{sub4sec4} how to sample from the posterior distributions of $\Delta
, R$ and $\psi$.
We use the parameter values generated to compute posterior
probabilities such as
$P\big(\Delta>\Delta_0\,|\,D\big), P\big(R>R_0\,|\,D\big)$, $P\big(\psi
>\psi_0\,|\,D\big)$,
and Bayesian credible\vadjust{\eject} intervals. The Bayes factors are introduced in
Section \ref{sub5sec4} to
compare the models under the hypotheses (i) $H_0$: $\lambda_0=\lambda
_1$ and $H_1$: $\lambda_1\neq\lambda_0$ and
(ii) $H^\ast_0$: $\gamma_1=\gamma_0$ versus $H^\ast_1$: $\gamma_1\neq
\gamma_0$.

\subsection{The Uniform Prior Distribution}\label{sub1sec4}

The prior $\pi_U(\gamma,\lambda_0,\lambda_1) = 2$ for $(\gamma,\lambda
_0,\lambda_1)\in\Omega$
refers to the uniform distribution under Dallal's model. This prior is
proper and expresses a complete
indifference of one vector of parameter values over another. The
resulting posterior distribution is
\begin{eqnarray}\label{eq1sub1sec4}
\pi_U\big(\gamma, \lambda_0, \lambda_1|D\big)
& \propto& 2^{m_{1+}} \gamma^{m_{1+}}
(1-\gamma)^{m_{2+}}\lambda_0^{m_{10} + m_{20}}
\big[1-(1+\gamma)\lambda_0\big]^{m_{00}}\lambda_1^{m_{11} + m_{21}}
\nonumber\\
& & \times\big[1-(1+\gamma)\lambda_1\big]^{m_{01}},
\qquad\qquad(\gamma,\lambda_0,\lambda_1)\in\Omega,\qquad\qquad
\end{eqnarray}
where $m_{1+}=m_{10}+ m_{11}$ and $m_{2+}=m_{20}+ m_{21}$.
As a result, the marginal posterior distribution of the nuisance
parameter $\gamma$ is
\begin{eqnarray}\label{eq2sub1sec4}
\pi_U\big(\gamma|D\big) = \dfrac{2^{m_{1+}+1}}{\B(m_{1+}+1,m_{2+}+1)}
\dfrac{\gamma^{m_{1+}}
(1-\gamma)^{m_{2+}}}{(1+\gamma)^{m_{1+}+m_{2+}+2}}, \qquad0<\gamma<1,
\end{eqnarray}
and the conditional posterior distribution of $\lambda_i, i=0,1$ given
$\gamma$ is
\begin{eqnarray}\label{eq3sub1sec4}
\pi_U\big(\lambda_i|\,\gamma, D\big)
& = & (1+\gamma)^{m_{1i}+m_{2i}+1}\dfrac{\lambda_i^{m_{1i}+m_{2i}}
\big[1-(1+\gamma)\lambda_i\big]^{m_{0i}}}{\B(m_{1i}+m_{2i}+1,m_{0i}+1)},
\quad0<\lambda_i<\dfrac{1}{1+\gamma},\qquad
\end{eqnarray}
where $\B(.,.)$ refers to the Beta function.
In other words, $\dfrac{1-\gamma}{1+\gamma}\sim\Be(m_{2+}+1,m_{1+}+1)$ and
$(1+\gamma)\lambda_i|\gamma\sim\Be(m_{1i}+m_{2i}+1,m_{0i}+1),\,i=0,1$,
where the notation
$\Be(\alpha,\beta)$ represents the standard Beta distribution with
shape parameters $\alpha$ and $\beta$.
The uniform prior can be viewed as a process of adding $1/2$ to the
summary statistics in the bottom four cells of the $3\times2$ table
and 1 to the top two cells. The uniform prior is appealing in
situations where the physical system imposes a natural parameterization
with a nice physical interpretation. In general, the uniform
distribution lacks the property of parameterization invariance.

\subsection{Jeffreys' Prior}\label{sub2sec4}

Jeffreys' prior has the property of being invariant under
a one-to-one reparameterization \citep{Jeffreys1946}. %(Jeffreys, 1946).
Regardless of the
parameterization used, Jeffreys' prior distribution is proportional to
the square
root of the absolute value of the determinant of the Fisher's
information matrix.

Define $U = (1+\gamma)\lambda_0$ or equivalently $\lambda_0 = \dfrac
{U}{1+\gamma}$,
and let $V = (1+\gamma)\lambda_1$ or equivalently $\lambda_1 =\dfrac
{V}{1+\gamma}$.
Under this new parametrization, the parameter space reduces to the
interval $(0,1)$
for each of the three parameters $\gamma$, $U$ and $V$. Thus, we are no
longer dealing
with a constrained parameter space. Moreover, the triplet $(\gamma, U,
V)$ forms\vadjust{\eject} a set of orthogonal
parameters in the sense of \citet{CoxReid1987}, %Cox and Reid (1987),
that is, the off-diagonal elements of the
expected Fisher information matrix are all zero. Propositions \ref
{prop1sub2sec4} and
\ref{prop2sub2sec4} give Jeffreys' priors for $(\gamma, U, V)$ and
$(\gamma, \lambda_0, \lambda_1)$, which are derived in Appendix 1. %
\begin{prop}\label{prop1sub2sec4}
Jeffreys' prior under the parameterization $(\gamma,U,V)$ is
\begin{equation}\label{eq1sub2sec4}
\pi_J(\gamma, u, v) \propto \dfrac{\gamma^{1/2-1}(1-\gamma
)^{1/2-1}}{(1+\gamma)} (u + rv)^{1/2}
u^{1/2-1}(1-u)^{1/2-1} v^{1/2-1}(1-v)^{1/2-1}
\end{equation}
for $0<\gamma, u, v<1$ and it is proper, where $r=\dfrac
{m_{+1}}{m_{+0}}$ is the ratio of the sample sizes in the
two treatment groups.
\end{prop}
Under Jeffreys' prior, the nuisance parameter, $\gamma$, is independent
of both $U$ and $V$ and its
marginal prior distribution is given by
\begin{eqnarray}\label{eq2sub2sec4}
\pi_J(\gamma) = \dfrac{\sqrt{2}}{\pi} \,\frac{\gamma^{1/2-1}
\;(1-\gamma)^{1/2-1}}{(1+\gamma)}, \qquad0<\gamma<1.
\end{eqnarray}
However, Jeffreys' prior depends indirectly on the sample sizes in both
groups through their ratio, $r$.

\begin{prop}\label{prop2sub2sec4}
In the original space, Jeffreys' prior reduces to
\begin{eqnarray}\label{eq3sub2sec4}
\pi_J(\gamma, \lambda_0, \lambda_1) & \propto&
\sqrt{\frac{(1+\gamma) (\lambda_0+r\lambda_1)}{\gamma
(1-\gamma) \lambda_0 \lambda_1
\big[1-(1+\gamma)\lambda_0\big] \big[1-(1+\gamma)\lambda_1\big]}},\\
& & (\gamma, \lambda_0, \lambda_1) \in\Omega.\nonumber
\end{eqnarray}
\end{prop}
The posterior distribution resulting from the use of Jeffreys' prior
is
\begin{align}\label{eq4sub2sec4}
\pi_J\big(\gamma, \lambda_0, \lambda_1|D\big) \propto&
\dfrac{2^{m_{1+}+1/2}}{\B(m_{1+}+1/2,m_{2+}+1/2)} \gamma^{m_{1+}+1/2-1}
(1-\gamma)^{m_{2+}+1/2-1}(1+\gamma)^{1/2}\nonumber\\
& \times(\lambda_0 + r\lambda_1)^{1/2}\,\lambda_0^{m_{10}+m_{20}+1/2-1}
\big[1-(1+\gamma)\lambda_0\big]^{m_{00}+1/2-1} \nonumber\\
& \times\lambda_1^{m_{11}+m_{21}+1/2-1}\big[1-(1+\gamma)\lambda_1\big
]^{m_{01}+1/2-1},
\qquad(\gamma,\lambda_0,\lambda_1)\in\Omega.
\end{align}

\subsection{Reference Priors}\label{sub3sec4}

Despite its success in the one-parameter context, Jeffreys'
non-informative prior methodology often
runs into serious difficulties in multiparameter problems \citep
{DattaGhosh1996}. %(Datta and Ghosh, 1996).
The prior
distribution may be difficult to derive and too complex to be easily
interpretable. This is often the case when
several nuisance parameters are present. For our problem, Jeffreys'
prior is difficult to interpret given that it
depends on the ratio of sample sizes $m_{+1}$ and $m_{+0}$. In this
context, the reference prior
may be more preferred \citep{Bernardo1979,BergerBernardo1989},
%%(Bernardo, 1979, Berger and Bernardo, 1989),
loosely defined as vague priors\vadjust{\eject} with the least amount of information.
Here, we divide the vector of parameters into two sets: parameters of
interest and nuisance parameters. Then, we consider the parameters
sequentially in the process
of deriving a reference prior. \citet{BergerBernardo1992a,BergerBernardo1992b,BergerBernardo1992c}
%Berger and Bernardo
%(1992a,b,c)
took this idea to another level by
suggesting to split the parameter vector into multiple groups according
to their orders of inferential
importance. The two authors provided a general algorithm for the
construction of reference priors. Hence,
reference priors are not uniquely defined. We derive the reference
prior for the parameters
$\gamma$, $U$, and $V$, leading to an induced reference prior for
$\gamma$, $\lambda_0$, and $\lambda_1$.

We start with the two group orderings: (i) $\{U,V\}$ and then $\{\gamma
\}$ and (ii)
$\{V,U\}$ and then $\{\gamma\}$.
\begin{prop}\label{prop1sub3sec4}
Bernardo's reference prior corresponding to the groups ordering $\{U,V\}
$ and then $\{\gamma\}$ or
$\{V,U\}$ and then $\{\gamma\}$ is
\begin{equation}\label{eq1sub3sec4}
\pi_R(\gamma, u, v) = \dfrac{2^{1/2}\gamma^{1/2-1}
(1-\gamma)^{1/2-1}}{\B(1/2,1/2)(1+\gamma)}
\dfrac{u^{1/2-1}(1-u)^{1/2-1}}{\B(1/2,1/2)} \dfrac
{v^{1/2-1}(1-v)^{1/2-1}}{\B(1/2,1/2)},\quad
\end{equation}
$0<\gamma, u, v<1$. That is, $\dfrac{1-\gamma}{1+\gamma}\sim\Be
(1/2,1/2)$,\,
$U\sim\Be(1/2,1/2)$ \,and\, $V\sim\break\Be(1/2,1/2)$.
\end{prop}
\begin{prop}\label{prop2sub3sec4}
In the original parameterization, $(\gamma, \lambda_0, \lambda_1)$,
Bernardo's reference
prior is equivalent to
\begin{align}\label{eq2sub3sec4}
\pi_R(\gamma, \lambda_0, \lambda_1) = & \;
\dfrac{\sqrt{2}}{\pi} \gamma^{1/2-1}
(1-\gamma)^{1/2-1} \dfrac{\lambda_0^{1/2-1}
\big[1-(1+\gamma)\lambda_0\big]^{1/2-1}}{\B(1/2,1/2)}\nonumber\\
& \; \times\dfrac{\lambda_1^{1/2-1}
\big[1-(1+\gamma)\lambda_1\big]^{1/2-1}}{\B(1/2,1/2)},\qquad
(\gamma, \lambda_0, \lambda_1) \ \in\Omega.
\end{align}
\end{prop}
The proofs of Propositions \ref{prop1sub3sec4} and \ref{prop2sub3sec4}
are given in Appendix 2. %\ref{sec1app}.

Unlike Jeffreys' prior, $\gamma$, $U$, and $V$ are independent under
Bernardo's reference prior.
The posterior distribution resulting from the use of Bernardo's
reference prior is
\begin{align*}
& \; \pi_R\big(\gamma, \lambda_0, \lambda_1|D\big) \\
\propto& \;
\dfrac{2^{m_{1+}+1/2}}{\B(m_{1+}+1/2,m_{2+}+1/2)}\gamma^{m_{1+}+1/2-1}
(1-\gamma)^{m_{2+}+1/2-1}\lambda_0^{m_{10}+m_{20}+1/2-1}\\
& \; \times\big[1-(1+\gamma)\lambda_0\big]^{m_{00}+1/2-1}
\lambda_1^{m_{11}+m_{21}+1/2-1}\big[1-(1+\gamma)\lambda_1\big]^{m_{01}+1/2-1},
\end{align*}
$(\gamma,\lambda_0,\lambda_1)\in\Omega$ or equivalently
\begin{align*}
& \; \pi_R\big(\gamma, u, v \,|D) \\
= & \; \dfrac{2^{m_{1+}+1/2}}{\B(m_{1+}+1/2, m_{2+}+1/2)}
\dfrac{\gamma^{m_{1+}+1/2-1} (1-\gamma)^{m_{2+}+1/2-1}}{(1+\gamma
)^{m_{1+}+m_{2+}+1}} \nonumber\\
& \; \times\dfrac{u^{m_{10}+m_{20}+1/2-1}(1-u)^{m_{00}+1/2-1}}{\B
(m_{10}+m_{20}+1/2, m_{00}+1/2)}
\dfrac{v^{m_{11}+m_{21}+1/2-1} (1-v)^{m_{01}+1/2-1}}{\B
(m_{11}+m_{21}+1/2, m_{01}+1/2)},
\end{align*}
$0<\gamma,u,v<1$. The reference prior can be viewed as adding $1/4$ to
each of the bottom four cells
of the $3\times2$ table and $1/2$ to the top two cells.

\citet{GhoshMukerjee1992} %Ghosh and Mukerjee (1992)
advise reversing the role of parameters of interest and
nuisance parameters to obtain a reverse reference prior. That is,
reconsider the group ordering
of $\{\gamma\}$ and then $\{U,V\}$ or $\{\gamma\}$ and then $\{V,U\}$.
The reference prior remains unchanged as shown in Appendix 2. %
We also discuss the idea of reference priors under the partial
information introduced by
\citet{SunBerger1998} %Sun and Berger (1998)
and known as conditional reference priors in Appendix 2. %\ref{sec1app}.

\subsection{Sampling from the Posterior Distribution}\label{sub4sec4}

With analytical solutions difficult to derive or compute, we turn to
Monte-Carlo simulation methods.
We first discuss how to simulate $\gamma$ from the distribution
\begin{eqnarray*}
f(\gamma) & = &
\dfrac{2^{\mu}}{\B(\mu,\nu)}\dfrac{\gamma^{\mu-1} \;(1-\gamma)^{\nu-1}}{
(1+\gamma)^{\mu+\nu}}\,,\qquad0<\gamma<1,
\end{eqnarray*}
which includes the marginal posterior distributions, $\pi_U(\gamma\,|\,
D)=\pi_J(\gamma\,|\,D)$, as special cases with
$\mu=m_{1+}+\tfrac{1}{2}$ and $\nu=m_{2+}+\tfrac{1}{2}$.
Then, we discuss how to simulate $(\lambda_1,\lambda_0)$ jointly from
$f(\lambda_1,\lambda_0 \mid D)$.

We propose two direct and efficient approaches to generate $\gamma$
from $f(\gamma)$.
\begin{itemize}
\item[(i)] Let $\gamma= \dfrac{e^\Psi}{1+e^\Psi}$. We show in Appendix
3 %\ref{sec3app}
that
$\phi= \Psi-\log(2)$ has the same distribution as $\mbox{logit}(p) =
\log\left(
\dfrac{p}{1-p}\right)$, where $p\sim\Be(\mu,\nu)$. Thus, simulate
$p_i\sim\Be(\mu,\nu)$, compute $\phi_i=\mbox{logit}(p_i)-\log(2)$, and set
$\gamma_i=\dfrac{e^{\phi_i}}{1+e^{\phi_i}},\,i=1,\dots,M$.
\item[(ii)] Let $\gamma= \dfrac{1-\pi}{1+\pi}$. We show in Appendix 3 %
that
$\pi\sim\Be(\nu,\mu)$. Thus, simulate $\pi_i\sim\Be(\nu, \mu)$ and set
$\gamma_i = \dfrac{1-\pi_i}{1+\pi_i},\,i=1,\dots,M$.
\end{itemize}

We now focus on the joint marginal posterior distribution of $(U,V)$
obtained under
Jeffreys' prior
\begin{align*}
f(u, v\,|\,D) \propto& \; (u + rv)^{1/2}\,
\dfrac{u^{m_{10}+m_{20}-1/2}(1-u)^{m_{00}-1/2}}{\B(m_{10}+m_{20}+1/2,
m_{00}+1/2)}\,\\
& \times\dfrac{v^{m_{11} + m_{21}-1/2} (1-v)^{m_{01}-1/2}}{\B
(m_{11}+m_{21}+1/2, m_{01}+1/2)},
\end{align*}
which is independent of $\gamma$. To simulate $M$ observations
$(u_i,v_i), \; i=1,\cdots,M$,
we proceed as follows:
\begin{itemize}
\item[(a)] Simulate independent observations $(u_i, v_i), \; i=1,\cdots
,M$, with $u_i\sim
\Be\big(m_{10}+m_{20}+1/2, m_{00}+1/2\big)$ and $v_i\sim
\Be\big(m_{11}+m_{21}+1/2, m_{01}+1/2\big)$.
\item[(b)] Compute the weights $w_i=(u_i+rv_i)^{1/2}, \; i=1,\cdots,M$.
\item[(c)] Use the acceptance/rejection sampling method: Simulate $\xi
_i\sim U(0,1)$ and
accept the pair $(u_i,v_i)$ only if $\xi_i < w_i/(1+r)$.
\item[(d)] Or use the importance sampling method, where all the pairs
$(u_i,v_i)$ are
accepted and use the weights $w_i$ to correct for the bias in the
computation of
posterior mean and quantiles.
\end{itemize}
Under the reference prior, we simulate independent observations $(u_i,
v_i), \; i=1,\cdots,M$, with $u_i\sim
\Be\big(m_{10}+m_{20}+1/2, m_{00}+1/2\big)$ and $v_i\sim
\Be\big(m_{11}+m_{21}+1/2, m_{01}+1/2\big)$. Having simulated a triplet
$(\gamma_i, U_i,V_i)$, we compute $\lambda^i_0=U_i/(1+\gamma_i)$ and
$\lambda^i_1=V_i/(1+\gamma_i)$ as well as the risk difference, $\Delta
_i=\dfrac{(V_i-U_i)}{(1+\gamma_i)}$,
the risk ratio, $R_i=\dfrac{V_i}{U_i}$ (does not depend on $\gamma$),
and the odds ratio, $\psi_i=\dfrac{V_i\big[1+\gamma_i-U_i\big]}{U_i\big
[1+\gamma_i-V_i\big]}$.
These simulated values are in turn used to compute posterior
probabilities and Bayesian credible intervals such as
equal-tailed intervals and highest posterior density (HPD) intervals.
Bayesian HPD intervals are the shortest intervals
containing the parameter of interest with the desired posterior
coverage probability. They are more desirable
than the commonly used equal-tailed intervals when the posterior
distribution is highly skewed,
but are more difficult to compute. \citet{ChenShao1999} %Chen and Shao
%(1999)
develop the Monte Carlo method to compute the HPD intervals.
In their paper, they also discussed how to compute Monte-Carlo-based
Bayesian credible intervals under importance sampling.

\subsection{Marginal Predictive Distribution and Bayes Factor}\label{sub5sec4}

In this section, we derive the marginal predictive distribution and the
Bayes factor under
Jeffreys' prior and Bernardo's prior. For simplicity and clarity of the
presentation,
we work with the parameterization $(\gamma, U, V)$. In addition, we
derive a
single Bayes factor formula for these two priors. To accomplish this
single formulation,
we use a family of prior distributions that encapsulates both Jeffreys'
prior and the reference prior as special cases.

\setcounter{secnumdepth}{3}
\subsubsection{$H_0$: $\lambda_0=\lambda_1$ versus $H_1$: $\lambda_1\neq
\lambda_0$}\label{sub1sub4sec4}

Consider the hypotheses $H_0$: $\lambda_0=\lambda_1=\lambda$ and $H_1$:
$\lambda_1\neq\lambda_0$ or equivalently
$H_0$: $U=V=\theta$ against $H_1$: $ U \neq V$. Under $H_1$, we
consider the family of prior distributions
\begin{align*}
& \; \pi_{H_1}\big(\gamma, u, v\big) \\
= & \; \dfrac{1}{K} \dfrac{2^{\frac{1}{2}}}{\B(\frac{1}{2},\frac{1}{2})}
\dfrac{\gamma^{\frac{1}{2}-1} (1-\gamma)^{\frac{1}{2}-1}}{(1+\gamma)} (u+rv)^{d}
\dfrac{u^{\frac{1}{2}-1}(1-u)^{\frac{1}{2}-1}}{\B(\frac{1}{2}, \frac
{1}{2})} \dfrac{v^{\frac{1}{2}-1}
(1-v)^{\frac{1}{2}-1}}{\B(\frac{1}{2}, \frac{1}{2})},
\end{align*}
where $K$ is the normalizing constant and $0<\gamma, u, v<1$. Note that
when $d=0$, $K=1$. Two choices of $d$
are of interest: $d=0$ corresponding to the reference prior and $d=1/2$
corresponding to Jeffreys' prior.
The marginal predictive distribution under $H_1$ is
\begin{align*}
& \; p_{H_1}(m_{10}, m_{20}, m_{11}, m_{21}) \\
= & \; \dfrac{1}{K} p(m_{10}, m_{20}, m_{11}, m_{21})
\int_0^1 \int_0^1 (u+rv)^{d}
\dfrac{u^{m_{10}+m_{20}+\frac{1}{2}-1} (1-u)^{m_{00}+\frac{1}{2}-1}}{\B
\big(m_{10}+m_{20}+\frac{1}{2},m_{00}+\frac{1}{2}\big)} \\
& \; \times
\dfrac{v^{m_{11}+m_{21}+\frac{1}{2}-1} (1-v)^{m_{01}+\frac{1}{2}-1}}{\B
\big(m_{11}+m_{21}+\frac{1}{2},m_{01}+\frac{1}{2}\big)} dv\,du,
\end{align*}
where
\begin{align*}
p(m_{10}, m_{20}, & m_{11}, m_{21}) = {m_{+0} \choose{m_{00}, m_{10}, m_{20}}}
{m_{+1} \choose{m_{01}, m_{11}, m_{21}}} \dfrac{\B\big(m_{1+}+\frac
{1}{2}, m_{2+}+\frac{1}{2}\big)}{\B\big(\frac{1}{2},\frac{1}{2}\big)} \\
& \; \times\dfrac{\B\big(m_{10}+m_{20}+\frac{1}{2}, m_{00}+\frac
{1}{2}\big)}{\B\big(\frac{1}{2},\frac{1}{2}\big)}
\dfrac{\B\big(m_{11}+m_{21}+\frac{1}{2}, m_{01}+\frac{1}{2}\big)}{\B\big
(\frac{1}{2},\frac{1}{2}\big)}.
\end{align*}
Under $H_0$, the likelihood reduces to
\[
L(\gamma, \theta) = \left\{\prod_{i=0}^1 {m_{+i} \choose{m_{0i},
m_{1i}, m_{2i}}}\right\} \;2^{m_{1+}} \dfrac{\gamma^{m_{1+}}
(1-\gamma)^{m_{2+}}}{(1+\gamma)^{m_{1+}+m_{2+}}}\,
\theta^{m_{1+}+m_{2+}}(1-\theta)^{m_{0+}},
\]
where $\theta=(1+\gamma)\lambda$. Jeffreys' and Bernardo's priors
belong to the family of priors
\begin{eqnarray*}
\pi_{H_0}(\gamma, \theta) & = & \dfrac{2^{\frac{1}{2}}}{\B(\frac
{1}{2},\frac{1}{2})}\,
\dfrac{\gamma^{\frac{1}{2}-1}\;(1-\gamma)^{\frac{1}{2}-1}}{(1+\gamma)}
\;
\dfrac{\theta^{a-1}\;(1-\theta)^{\frac{1}{2}-1}}{\B(a,\frac
{1}{2})},\quad0<\gamma,\theta<1.
\end{eqnarray*}
%
%$a=1$ is used in Jeffreys' prior and $a=1/2$ is used in the reference
%prior.
Thus, the resulting marginal predictive distribution is
\begin{align*}
& \; p_{H_0}(m_{10}, m_{20}, m_{11}, m_{21}) \\
= & \; {m_{+0} \choose{m_{00},
m_{10}, m_{20}}}\;{m_{+1} \choose{m_{01},m_{11}, m_{21}}}\,
\dfrac{\B\big(m_{1+}+\frac{1}{2},\, m_{2+}+\frac{1}{2}\big)}{\B\big
(\frac{1}{2},\frac{1}{2}\big)} \\
& \times\dfrac{\B\big(m_{2+}+m_{1+}+a,\,m_{0+}+\frac{1}{2})}{\B(a,\frac
{1}{2}\big)}.
\end{align*}
The ratio of these two marginal predictive distributions, $BF^\lambda
_{01}$, under the condition
$P(H_1)=P(H_0)=\tfrac{1}{2}$ (the Bayes factor for testing $H_0$ vs
$H_1$) satisfies
\begin{align*}
& \dfrac{1}{BF_\lambda}
= \dfrac{\B\big(m_{10}+m_{20}+\frac{1}{2}, m_{00}+\frac{1}{2}\big) \B
\big(m_{11}+m_{21}+\frac{1}{2}, m_{01}+\frac{1}{2}\big) \B\big(a,\frac
{1}{2}\big)}
{\B\big(m_{2+}+m_{1+}+a, m_{0+}+\frac{1}{2}\big) \B\big(\frac
{1}{2},\frac{1}{2}\big) \B\big(\frac{1}{2},\frac{1}{2}\big)}\\
& \times\int_0^1 \int_0^1 \dfrac{(u+rv)^{d}}{K} \dfrac
{u^{m_{10}+m_{20}+\frac{1}{2}-1}
(1-u)^{m_{00}+\frac{1}{2}-1}}{\B(m_{10}+m_{20}+\frac{1}{2},m_{00}+\frac
{1}{2})} \dfrac{v^{m_{11}+m_{21}+\frac{1}{2}-1}
(1-v)^{m_{01}+\frac{1}{2}-1}}{\B(m_{11}+m_{21}+\frac{1}{2},m_{01}+\frac
{1}{2})} dv du.
\end{align*}
Under Jeffreys' prior, the constant $K$ and the integral term in the
Bayes factor are computed using computer simulation.
Under the reference prior, the integral term disappears and the Bayes
factor is computed exactly using only the
Beta functions.

\subsubsection{$H^\ast_0$: $\gamma_1=\gamma_0$ versus
$H^\ast_1$: $\gamma_1\neq\gamma_0$}\label{sub2sub4sec4}

One of the statements made in Dallal's model is that the parameter
$\gamma$ is constant.
As discussed earlier, this assumption can be relaxed to $P\big
(Z_{ijk}=1 \,|\, Z_{ij(3-k)}=1\big) = 1-\gamma_i$, giving rise to
the full or saturated model. Therefore, it is important to test the
hypothesis $H^\ast_0$: $\gamma_1=\gamma_0=\gamma$
(Dallal's reduced model) versus the alternative hypothesis $H^\ast_1$:
$\gamma_1\neq\gamma_0$ (Dallal's full model). Under $H^\ast_1$, $U$ and $V$
are redefined as follows: $U=(1+\gamma_0)\lambda_0$ and $V=(1+\gamma
_1)\lambda_1$.
Under $H^\ast_0$, the prior is $\pi_{H_1}(\gamma,u,v)$ defined in
Section \ref{sub1sub4sec4}.
Under $H^\ast_1$, the family of priors under consideration is
\begin{align*}
\pi_{H^\ast_1}(\gamma_0, \gamma_1, u, v) = &
\dfrac{2^{\frac{1}{2}}\gamma_0^{\frac{1}{2}-1}\,
(1-\gamma_0)^{\frac{1}{2}-1}}{\B(\frac{1}{2},\frac{1}{2})(1+\gamma_0)}\;
\dfrac{2^{\frac{1}{2}}\gamma_1^{\frac{1}{2}-1}\,
(1-\gamma_1)^{\frac{1}{2}-1}}{\B(\frac{1}{2},\frac{1}{2})(1+\gamma_1)}\,
\dfrac{u^{a_0-1}\,(1-u)^{\frac{1}{2}-1}}{\B(a_0,\frac{1}{2})}\\
&\, \times\dfrac{v^{a_1-1}\,(1-v)^{\frac{1}{2}-1}}{\B(a_1,\frac
{1}{2})},\,\qquad0<\gamma_0, \gamma_1, u, v <1.
\end{align*}
The resulting marginal predictive distribution is
\begin{align*}
p_{H^\ast_1}\big(m_{10}, m_{20}, m_{11}, m_{21}\big) = & {m_{+0}
\choose{m_{00},
m_{10}, m_{20}}}\;{m_{+1} \choose{m_{01}, m_{11}, m_{21}}}\,\dfrac{\B
\big(m_{10}+\frac{1}{2},
m_{20}+\frac{1}{2}\big)}{\B\big(\frac{1}{2},\frac{1}{2}\big)} \\
& \times\dfrac{\B\big(m_{11}+\frac{1}{2}, m_{21}+\frac{1}{2}\big)}{\B
\big(\frac{1}{2},\frac{1}{2}\big)}\;
\dfrac{\B(m_{10}+m_{20}+a_0, m_{00}+\frac{1}{2})}{\B\big(a_0,\frac
{1}{2}\big)}\\
& \times\dfrac{\B(m_{11}+m_{21}+a_1, m_{01}+\frac{1}{2})}{\B\big
(a_1,\frac{1}{2}\big)}.
\end{align*}
We are only concerned with two sets of parameters choices:
$d=1/2, a_0=a_1=1$ used in Jeffreys' prior and $d=0, a_0=a_1=1/2$ used
in the reference prior.

The Bayes factor for testing $H^\ast_0$ vs $H^\ast_1$ under the
condition $P(H_1)=P(H_0)=\tfrac{1}{2}$ is
\begin{align*}
& BF_\gamma = \dfrac{p_{H_1}\big(m_{10}, m_{20}, m_{11}, m_{21}\big
)}{p_{H^\ast_1}
\big(m_{10}, m_{20}, m_{11}, m_{21}\big)} \\
= & \; \dfrac{\B\big(m_{1+}+\frac{1}{2}, m_{2+}+\frac{1}{2}\big) \B\big
(a_0,\frac{1}{2}\big)\,\B\big(a_1,\frac{1}{2}\big)}
{\B\big(m_{10}+\frac{1}{2}, m_{20}+\frac{1}{2}\big) \B\big(m_{11}+\frac
{1}{2}, m_{21}+\frac{1}{2}\big)
\B\big(\frac{1}{2},\frac{1}{2}\big)}\\
& \times\dfrac{\B\big(m_{10}+m_{20}+\frac{1}{2}, m_{00}+\frac{1}{2}\big)
\B\big(m_{11}+m_{21}+\frac{1}{2}, m_{01}+\frac{1}{2}\big)}
{\B\big(m_{10}+m_{20}+a_0, m_{00}+\frac{1}{2}\big)
\B\big(m_{11}+m_{21}+a_1, m_{01}+\frac{1}{2}\big)}\\
& \times\int_0^1 \int_0^1 \Big\{ \dfrac{(u+rv)^{d}}{K} \dfrac
{u^{m_{10}+m_{20}+\frac{1}{2}-1}
(1-u)^{m_{00}+\frac{1}{2}-1}}{\B\big(m_{10}+m_{20}+\frac
{1}{2},m_{00}+\frac{1}{2}\big)}\\
& \hspace*{0.3in} \times\dfrac{v^{m_{11}+m_{21}+\frac{1}{2}-1}
(1-v)^{m_{01}+\frac{1}{2}-1}}{\B\big(m_{11}+m_{21}+\frac
{1}{2},m_{01}+\frac{1}{2}\big)} \Big\} dvdu.
\end{align*}
Under Jeffreys' prior, $BF_\gamma$ is computed using computer
simulation while under the reference prior it is
computed exactly.

%%----------------------------------------------------------------------------------------------
% Section 5
%%----------------------------------------------------------------------------------------------

\section{Comparisons of Bayesian and Frequentist Intervals: An
Empirical Study}\label{sec5}

%
% Figure 1: 4 sets of graphs
%

% \epsfig{figure=MSEPlot1.eps, height=6.0in, width=5.0in}
%the case $\Delta=0$.}

% \epsfig{figure=MSEPlot2.eps, height=6.0in, width=5.0in}
%from case $\Delta=0.8$.}

% \epsfig{figure=LengthPlot1.eps, height=6.0in, width=5.0in}
%generated from the case $\Delta=0.0$.}

% \epsfig{figure=LengthPlot2.eps, height=6.0in, width=5.0in}
%generated from the case $\Delta=0.0, 0.5, 0.8$.}

In this section, we investigate small, moderate and large-sample
performances of frequentist confidence
intervals (FCIs) and Bayesian credible intervals (BCIs) under three
criteria. For a set values for the model
parameters, 10,000 $3\times2$ bilateral data tables are generated from
the product of trinomial
distributions under a balanced design. The essence of these criteria
rely on the
principle that good FCIs (Wald FCIs described in Appendix 1 of the
Supplementary Web Materials) or good HPD
BCIs should have their true coverage close to or preferably larger than
the nominal value.
Indeed, FCIs and BCIs with deflated true coverage are recommended
against. Lengths of the intervals must
also be considered. We use the three following criteria:
\begin{itemize}
\item[(i)] the expected true coverage probability (ETCP) of the
interval $\big(\widehat{\Delta}_L,
\widehat{\Delta}_U\big)$ for $\Delta$, $P\big(\widehat{\Delta}_L\leq
\Delta\leq\widehat{\Delta}_U\big)$,
in repeated sampling;
\item[(ii)] the expected width (EWCI) of the interval of $\big(\widehat
{\Delta}_L, \widehat{\Delta}_U\big)$ for $\Delta$,
$\widehat{\Delta}_U-\widehat{\Delta}_L$, in repeated sampling; and
\item[(iii)] the expected mean square error (MSE) that reflects a
compromise between bias and precision, in repeated sampling.
\end{itemize}
Although these criteria are defined in terms of $\Delta$, we also
examine the behaviors of
these criteria for $(\gamma,\lambda_0,\lambda_1)$ as well. For the case
of the MSE, we go further by adding the three MSEs
corresponding to $(\gamma,\lambda_0,\lambda_1)$ to obtain a global
measure of MSE. Four equal sample size scenarios are chosen:
$m_{+0}=m_{+1}=m=10, 25, 50, 100$ to reflect small, moderate and large
sample size situations. For symmetry,
we only consider cases where $\Delta$ is non-negative. More
specifically, we examine cases with
$\Delta=\Delta_h=h/10, \,h \in\{ 0, 1, 2, 3, 4, 5, 6, 7, 8, 9\}$. Each
choice of $\Delta$ implies the constraint:
$0<\gamma<\gamma_{h,\max}=\min(1,1/\Delta_h-1)$. We then choose a grid
of $\gamma$ points:
$\gamma_{hj}=j\gamma_{h, \max}/10$ with $j=1,2,\dots,9$. These values
are meant to capture a wide range of behaviours of the
conditional probability of an occurrence of a particular characteristic
at one site given an occurrence of that
characteristic at the other site in the range $(1-\gamma_{h,\max},1)$.
Then, set $\lambda_{hj, \max}=\dfrac{1}{1+\gamma_{hj}}-\Delta_h$.
For the pair $(\lambda_0,\lambda_1)$, we select the values: $\lambda
_{0hjk} = k\lambda_{hj,\max}/10$ with $k=1,2,\dots,9$ and $\lambda
_{1hjk} =
\lambda_{0hjk}+\Delta_h$. Hence, for each sample size and $\Delta_h$,
we compute the three criterion functions for 81
combinations of the triplets $(\gamma, \lambda_0, \lambda_1)$. In the
Bayesian framework, we examine HPD
intervals under the uniform, Jeffreys', and reference priors.

The probability of obtaining degenerate results in the frequentist
framework (non-estimable model parameters
or model parameters on the boundary of the parameter space or
parameters with a zero variance)
is relatively high in smaller samples, and even higher when combined
with large and small $\gamma$ values.
These cases are eliminated from the frequentist calculations. The five
tables in Appendix 2 of the Supplementary
Web Materials give a summary of results of the simulation study. In
each cell, we compute
two numbers: the proportion of empirical coverages within 0.01 of the
nominal coverage and
the proportion of empirical coverages above -0.02 of the nominal
coverage. The latter summary
carries out more value than the former one. The findings for the ETCP
criterion are summarized as follows.
\begin{itemize}
\item[(a)] For $\Delta=0$, Bayesian HPD intervals perform similarly for
moderate and large sample sizes ($m=50, 100$) and the
proportions of coverages above -0.02 of the nominal coverage are close
to 100\% for $\lambda_0,\lambda_1$ and $\Delta$.
This finding actually holds true for $\Delta\leq0.4$.
\item[(b)] Wald FCIs perform poorly for estimating $\lambda_0$ and
$\lambda_1$
when $m=10, 25$ and $\Delta\leq0.3$.
\item[(c)] Overall, in terms of coverage probability, the uniform
distribution seems to perform better than the other methods when $\Delta
\leq
0.5$. When $\Delta\geq0.6$, Jeffreys' prior and the reference prior
perform better when estimating
$\lambda_0$, $\lambda_1$ and $\Delta$.
\item[(d)] When it comes to estimating $\gamma$, the uniform prior
outperforms the other methods regardless of the nominal value of $\Delta$.
\item[(e)] In general, Jeffreys' prior and the reference prior tend to
give similar results regardless of the nominal value of $\Delta$.
\end{itemize}

\begin{figure}[t!]
\includegraphics{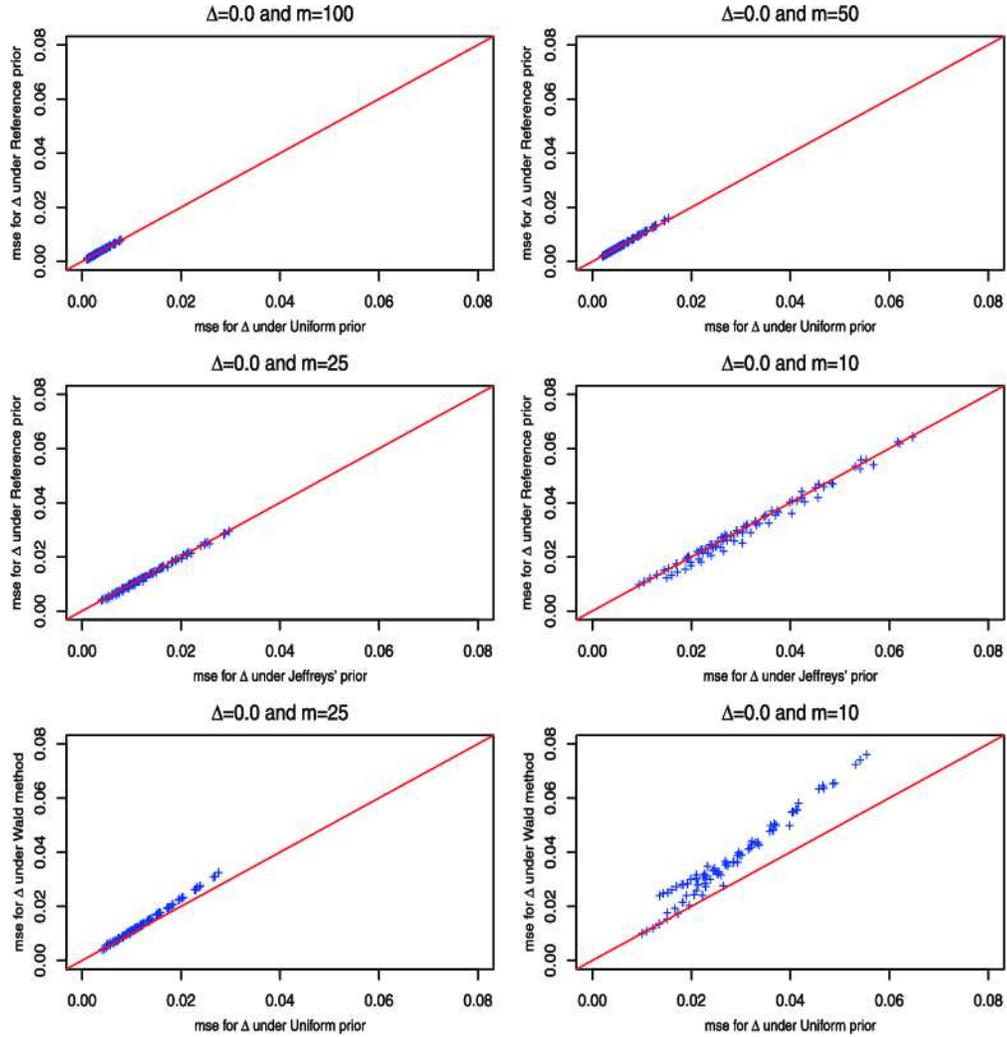}
\caption{Graphs of the 81 empirical MSEs for $\Delta$ generated from
the case $\Delta=0$.}
\end{figure}

\begin{figure}[t!]
\includegraphics{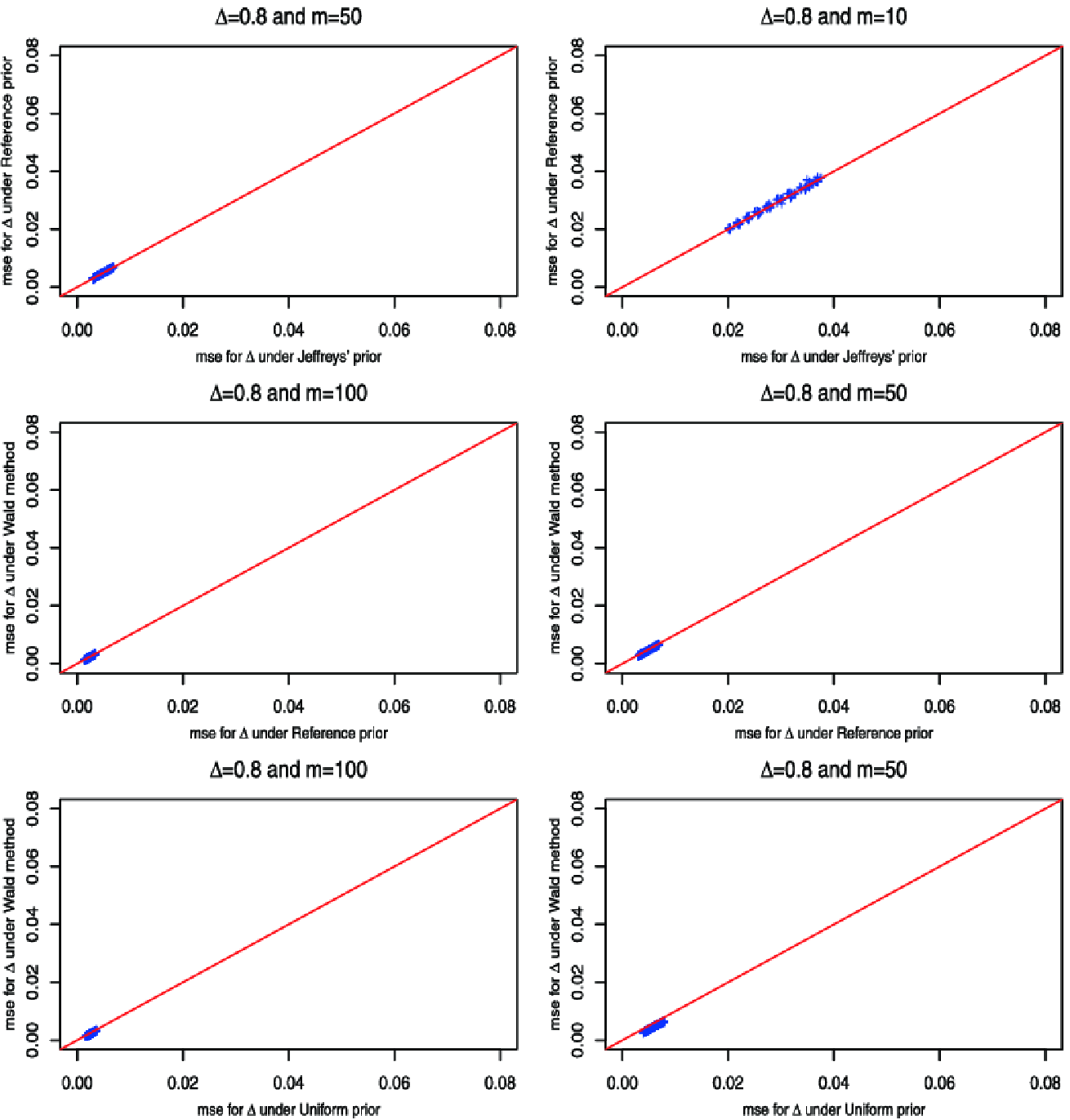}
\caption{Graphs of the 81 of empirical MSEs for $\Delta$ generated from
case $\Delta=0.8$.}
\end{figure}

\begin{figure}[t!]
\includegraphics{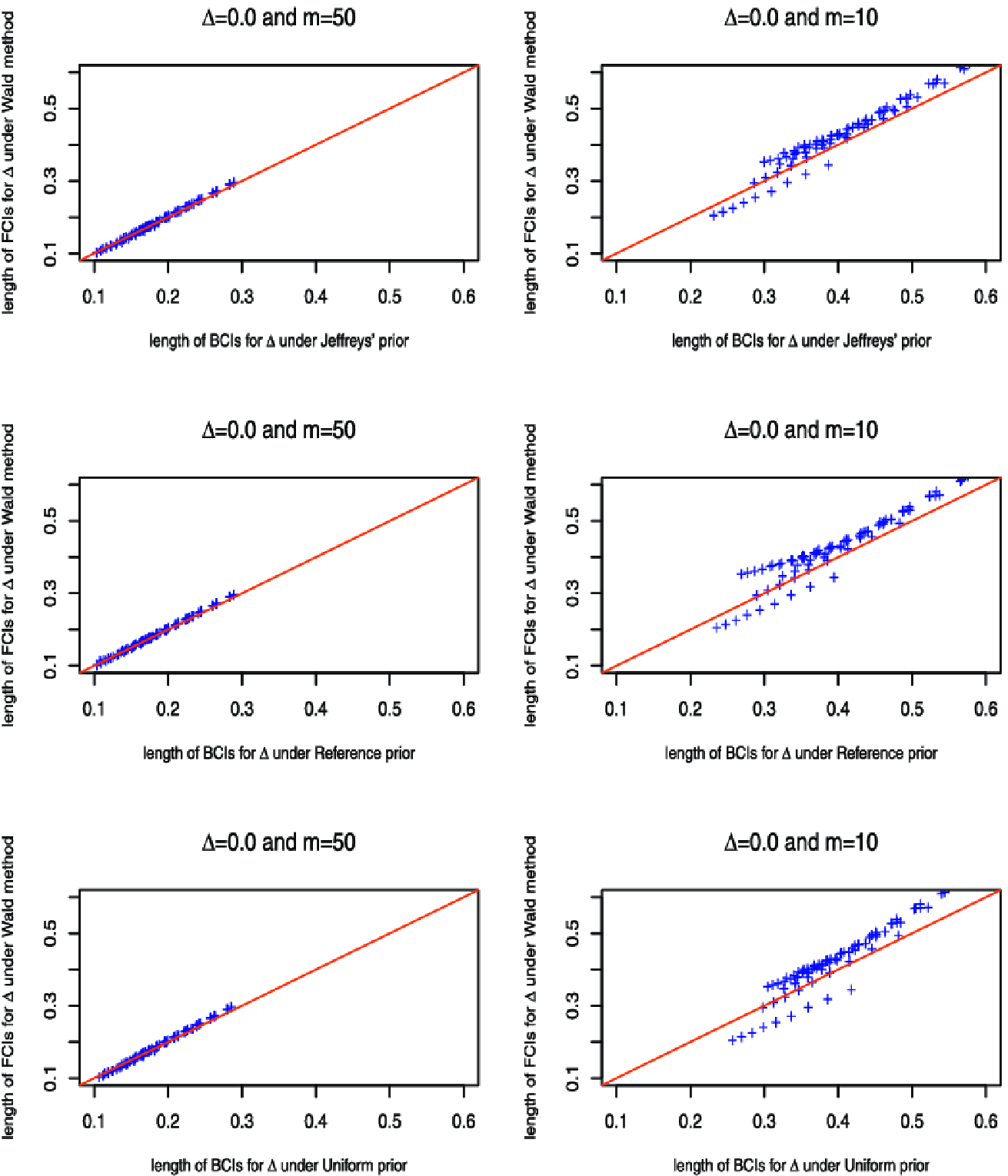}
\caption{Graph of the 81 empirical lengths of 90\% HPD intervals
generated from the case $\Delta=0.0$.}
\end{figure}

\begin{figure}[t!]
\includegraphics{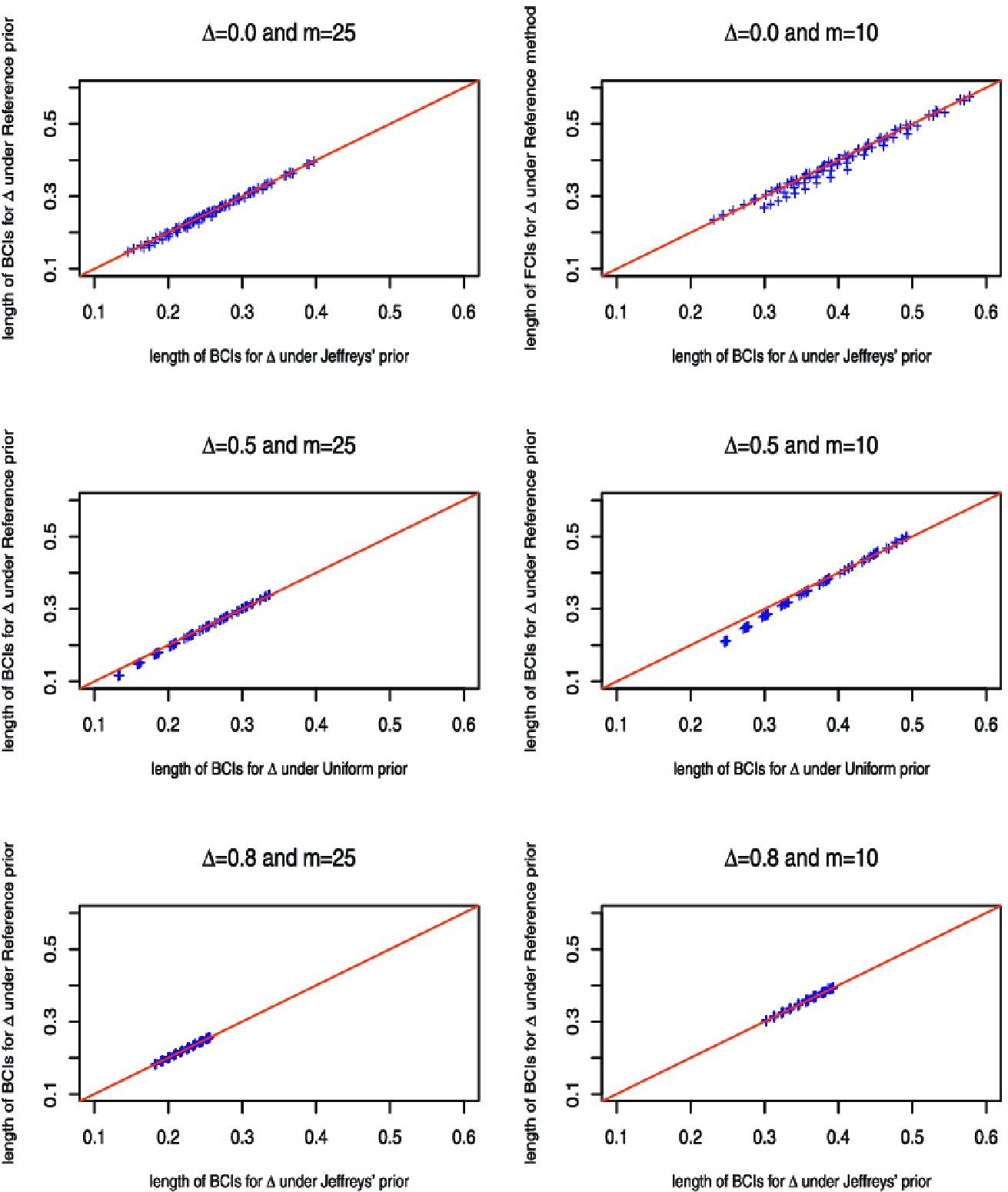}
\caption{Graph of the 81 empirical lengths of 90\% HPD intervals
generated from the case $\Delta=0.0, 0.5, 0.8$.}\vspace*{-3pt}
\end{figure}

A summary of the findings for the MSE criterion is given as follows.
\begin{itemize}
\item[(a)] Jeffreys' prior and the reference prior perform similarly
regardless of the nominal value of $\Delta$. The
MSEs for $\Delta$ are scattered around the line $y=x$, with increasing
deviations as $n$ gets smaller.
As $\Delta$ increases, the deviations from the 45 degree line also
increase with a slim advantage of Jeffreys' prior when
$\Delta>0.5$. See Figures 1 and 2. The same pattern is observed with
the global measure of MSE obtained
by adding the MSEs for $\gamma, \lambda_0$ and $\lambda_1$.
\item[(b)] Jeffreys' prior and the uniform prior have equivalent
properties for moderate and large sample
sizes when $\Delta\leq0.6$. However, when $\Delta\geq0.7$, the
superiority of Jeffreys' prior is highly evident.
\item[(c)] Wald FCIs perform worst when $\Delta\leq0.6$, with
differences worsening as $n$ gets smaller.
This behavior is reversed when $\Delta\geq0.7$. See Figures 1 and 2.
\item[(d)] In terms of the MSEs corresponding to $\gamma$, all the
Bayesian approaches perform better than
Wald's approach when $\Delta\leq0.8$.
\end{itemize}

The findings for the EWCI criterion are the following.
\begin{itemize}
\item[(a)] Jeffreys' prior and the reference prior give similar results
when estimating $\Delta$. See Figure 4.
The reference prior tends to be better than the uniform prior when
$\Delta\geq0.6$ and inferior when $\Delta\leq0.5$.
The Wald FCI is the worst when $\Delta\leq0.6$ and the best when
$\Delta\geq0.7$. See Figure 3.
\item[(b)] Wald FCIs perform worst when estimating $\gamma$. The
uniform prior is no better than Jeffreys' prior with
a slight edge to Jeffreys' prior when the sample size is small.
\end{itemize}

%%-------------------------------------------------------------------------------------------------------------------------------------------
% Section 6
%%-------------------------------------------------------------------------------------------------------------------------------------------

\section{Case Studies}\label{sec6}

\subsection{Bayesian Analysis of Bilateral Data with Sparse Data}\label
{sub1sec6}

%Mandel et al. (1982)
\citet{MandelBluestone1982}
considered a double-blind randomized clinical trial which compared the
antibiotics Cefaclor and Amoxicillin for the treatment of otitis media
with effusion (OME).
Among a total of 214 participants in the trial, 11 children were at
least six years old and
underwent bilateral tympanocentesis prior to randomization into one of
two groups (Cefaclor or
Amoxicillin). Children in each treatment group received a 14-day course
of treatment with one
of the antibiotics and dichotomous ear outcomes were determined (i.e.,
cured or not-cured)
and recorded. Table \ref{tab1sub1sec6} provides a summary of the data collected.
The primary goal of this investigation is to test if the cure rates
were identical between
Cefaclor and Amoxicillin. A further goal is to estimate the size of the
difference of the percentage change, $R-1$, in the
performance of the two medications. This example was discussed in
\citet{TangTangQiu2008} %Tang, Tang and Qiu (2008)
under Rosner's
model. They obtained $\widehat{\lambda}_0=0.875$ and $\widehat{\lambda
}_1=0.857$. They also found that
there is no evidence to reject the null hypothesis of equal cure rates.

\begin{table}[hbpt]
\caption{OME status after 14 days of Cefaclor or Amoxicillin treatments.}\label{tab1sub1sec6}
\smallskip\centering%\footnotesize
\begin{tabular}{|lllcccc|}\hline
{Number of ears with} &&& \multicolumn{4}{c|}{Treatment Group}\\\cline{3-7}
{OME being cured} &&& Amoxicillin &&& Cefaclor \\\hline
0 &&& 1 &&& 0\\
1 &&& 0 &&& 1\\
2 &&& 6 &&& 3\\
{\bf Total} &&& {\bf7} &&& {\bf4}\\\hline
\end{tabular}
\end{table}

Two striking features characterize this data. (i) The total sample
sizes in the Cefaclor and Amoxicillin groups are extremely small.
(ii) Table \ref{tab1sub1sec6} is sparse with a zero cell in each group.
As a result,
the normal approximation used in classical analysis does not apply. In
addition, the
maximum likelihood estimates (MLEs) of Dallal's reduced and saturated
models sit on the
boundary of the parameter space (leading to success probabilities not
truly allowed by the product
trinomial model). It is impossible to carry out frequentist inference
or apply a likelihood-based
model selection procedure. In this situation, a common ad-hoc\vadjust{\eject}
adjustment is to add
1/2 to each cell count. A problem with this ad-hoc adjustment is that
the total sample
sizes are not integer numbers. The normal approximation still does not
apply despite this
adjustment. So the Bayesian methodology appears here to be one of the
few alternatives.

\begin{table}[hbpt]
\caption{MLEs under Dallal's reduced and saturated models.}\label{tab4sub1sec6}
\smallskip
\centering%\footnotesize
\begin{tabular}{|c|c|c|c|c|c|c|c|}\hline
\multicolumn{8}{|c|}{{\bf Reduced Model}}\\\hline
\multicolumn{8}{|c|}{{\bf No Adjustments}}\\\hline
$U$ & $V$ & $\gamma$ & $\lambda_0$ & $\lambda_1$ & $\Delta$ & $R$ &
$\psi$ \\
1 & 6/7 & 1/19 & 0.95 & 57/70 & -19/140 & 6/7 & 3/13 \\\hline
\multicolumn{8}{|c|}{{\bf Ad-hoc Adjustment}}\\\hline
$U$ & $V$ & $\gamma$ & $\lambda_0$ & $\lambda_1$ & $\Delta$ & $R$ &
$\psi$ \\\hline
10/11 & 14/17 & 1/11 & 5/6 & 77/102 &-4/51 & 77/85 & 77/125\\\hline
\end{tabular}

\vspace{.1in}

\begin{tabular}{|c|c|c|c|c|c|c|c|c|}\hline
\multicolumn{9}{|c|}{{\bf Saturated Model}}\\\hline
\multicolumn{9}{|c|}{{\bf No Adjustments}}\\\hline
$U$ & $V$ & $\gamma_0$ & $\gamma_1$ & $\lambda_0$ & $\lambda_1$ &
$\Delta$ & $R$ & $\psi$\\\hline
1 & 6/7 & 1/7 & 0 & 0.875 & 6/7 & -1/56 & 48/49 & 6/7\\\hline
\multicolumn{9}{|c|}{{\bf Ad-hoc Adjustment}}\\\hline
$U$ & $V$ & $\gamma_0$ & $\gamma_1$ & $\lambda_0$ & $\lambda_1$ &
$\Delta$ & $R$ & $\psi$\\
10/11 & 14/17 & 3/17 & 1/27 & 17/22& 27/34 & 4/187 & 297/289 & 135/119\\
\hline
\end{tabular}
\end{table}

\begin{table}[hbpt]\caption{Posterior estimates and credible intervals
based on Dallal's model.}\label{tab2sub1sec6}
\smallskip\centering%\footnotesize
\begin{tabular}{|c|ccc|c|ccc|}\hline
\multicolumn{8}{|c|}{{\bf Jeffreys' Prior}} \\
\multicolumn{8}{|c|}{$P(\Delta<0|D)=0.739$, DIC$=$9.843, pD$=$1.613}
\\ \hline
% & \multicolumn{3}{|c|}{DIC=9.843, $pD=1.613$} \\\cline{2-4}
& mean & std & 95\% HPD & & mean & std & 95\% HPD \\\hline
$U$ & 0.904 & 0.119 & (0.649, 1.000) & $\lambda_1$ & 0.764 & 0.125 &
(0.515, 0.970) \\
$V$ & 0.819 & 0.127 & (0.570, 1.000) & $\Delta$ &-0.079 & 0.162 &
(-0.417, 0.265) \\
$\gamma$ & 0.076 & 0.061 & (0.000, 0.198) & $R$ & 0.931 & 0.245 &
(0.468, 1.370) \\
$\lambda_0$ & 0.842 & 0.120 & (0.597, 0.998) & $\psi$ & 1.013 & 2.052 &
(0.002, 3.423) \\ \hline
%$\lambda_1$ & 0.764 & 0.125 & (0.515, 0.970) \\
%$\Delta$ &-0.079 & 0.162 & (-0.417, 0.265)\\
%$R$ & 0.931 & 0.245 & (0.468, 1.370) \\
%$\psi$ & 1.013 & 2.052 & (0.002, 3.423) \\\hline
\multicolumn{8}{|c|}{{\bf Bernardo's Prior}} \\
\multicolumn{8}{|c|}{$P(\Delta<0|D)=0.737$, DIC$=$9.953, pD$=$1.631 }
\\ \hline
%& \multicolumn{3}{|c|}{DIC=9.953, $pD=1.631$}\\\cline{2-4}
& mean & std & 95\% HPD & & mean & std & 95\% HPD \\\hline
$U$ & 0.899 & 0.124 & (0.635, 1.000) & $\lambda_1$ & 0.756 & 0.128 &
(0.505, 0.969) \\
$V$ & 0.812 & 0.130 & (0.559, 0.999) & $\Delta$ & -0.082 & 0.168 &
(-0.439, 0.261) \\
$\gamma$ & 0.077 & 0.062 & (0.000, 0.199) & $R$ & 0.929 & 0.260 &
(0.469, 1.395) \\
$\lambda_0$ & 0.838 & 0.124 & (0.583, 0.997) & $\psi$ & 1.030 & 2.166 &
(0.001, 3.486) \\ \hline
%$\lambda_1$ & 0.756 & 0.128 & (0.505, 0.969) \\
%$\Delta$ & -0.082 & 0.168 & (-0.439, 0.261) \\
%$R$ & 0.929 & 0.260 & (0.469, 1.395) \\
%$\psi$ & 1.030 & 2.166 & (0.001, 3.486) \\\hline
\multicolumn{8}{|c|}{{\bf Uniform Prior}} \\
\multicolumn{8}{|c|}{$P(\Delta<0|D)=0.640$, DIC$=$10.378, pD$=$1.316}
\\ \hline
%& \multicolumn{3}{|c|}{DIC=10.378, $pD=1.316$} \\\cline{2-4}
& mean & std & 95\% HPD & & mean & std & 95\% HPD \\\hline
$U$ & 0.833 & 0.141 & (0.549, 1.000) & $\lambda_1$ & 0.713 & 0.127 &
(0.461, 0.935) \\
$V$ & 0.778 & 0.131 & (0.524, 0.989) & $\Delta$ & -0.051 & 0.177 &
(-0.400, 0.317) \\
$\gamma$ & 0.095 & 0.066 & (0.003, 0.224) & $R$ & 0.973 & 0.303 &
(0.453, 1.536) \\
$\lambda_0$ & 0.764 & 0.136 & (0.491, 0.977) & $\psi$ & 1.200 & 1.766 &
(0.016, 3.796) \\ \hline
%$\lambda_1$ & 0.713 & 0.127 & (0.461, 0.935) \\
%$\Delta$ & -0.051 & 0.177 & (-0.400, 0.317) \\
%$R$ & 0.973 & 0.303 & (0.453, 1.536) \\
%$\psi$ & 1.200 & 1.766 & (0.016, 3.796) \\\hline
\end{tabular}
\end{table}
\begin{table}[hbt]\caption{Posterior estimates and credible intervals
based on the saturated model.}\label{tab3sub1sec6}
\smallskip\centering%\footnotesize
\begin{tabular}{|c|ccc|c|ccc|}\hline
\multicolumn{8}{|c|}{{\bf Jeffreys' Prior}}\\
\multicolumn{8}{|c|}{$P(\Delta>0|D)=0.550$, $P(\gamma_1-\gamma
_0<0|D)=0.890$, DIC$=$8.680, pD$=$1.619 } \\ \hline
%& \multicolumn{3}{|c|}{$P(\gamma_1-\gamma_0<0|D)=0.890$}\\
%& \multicolumn{3}{|c|}{DIC=8.680, pD=1.619}\\\cline{2-4}
& mean & std & 95\% HPD & & mean & std & 95\% HPD \\\hline
$U$ & 0.910 & 0.112 & (0.670, 1.000) & $\Delta$ & 0.021 & 0.179 &
(-0.342, 0.371) \\
$V$ & 0.824 & 0.123 & (0.583, 0.998) & $R$ & 1.065 & 0.299 & (0.548,
1.628) \\
$\gamma_0$ & 0.192 & 0.145 & (0.000, 0.481) & $\psi$ & 2.633 & 7.094 &
(0.001, 8.849) \\
$\gamma_1$ & 0.039 & 0.055 & (0.000, 0.151) & $\delta$ & -0.132 & 0.167
& (-0.491, 0.200) \\
$\lambda_0$ & 0.773 & 0.128 & (0.532, 0.988) & $\Delta\gamma$ &-0.153 &
0.155 & (-0.495, 0.114) \\
$\lambda_1$ & 0.795 & 0.125 & (0.553, 0.993) & & & & \\ \hline
%$\Delta$ & 0.021 & 0.179 & (-0.342, 0.371)\\
%$R$ & 1.065 & 0.299 & (0.548, 1.628) \\
%$\psi$ & 2.633 & 7.094 & (0.001, 8.849) \\
%$\delta$ & -0.132 & 0.167 & (-0.491, 0.200) \\
%$\Delta\gamma$ &-0.153 & 0.155 & (-0.495, 0.114) \\\hline
\multicolumn{8}{|c|}{{\bf Bernardo's Prior}}\\
\multicolumn{8}{|c|}{$P(\Delta>0|D)=0.541$, $P(\gamma_1-\gamma
_0<0|D)=0.889$, DIC$=$8.858, pD$=$1.645 }\\ \hline
%& \multicolumn{3}{|c|}{$P(\gamma_1-\gamma_0<0|D)=0.889$} \\
%& \multicolumn{3}{|c|}{DIC=8.858, pD=1.645} \\\cline{2-4}
& mean & std & 95\% HPD & & mean & std & 95\% HPD \\\hline
$U$ & 0.900 & 0.122 & (0.638, 1.000) & $\Delta$ & 0.019 & 0.186 &
(-0.361, 0.393) \\
$V$ & 0.813 & 0.129 & (0.560, 0.999) & $R$ & 1.067 & 0.326 & (0.508,
1.667) \\
$\gamma_0$ & 0.193 & 0.145 & (0.000, 0.484) & $\psi$ & 2.574 & 6.512 &
(0.004, 8.772) \\
$\gamma_1$ & 0.040 & 0.055 & (0.000, 0.153) & $\delta$ & -0.135 & 0.176
& (-0.502, 0.221) \\
$\lambda_0$ & 0.765 & 0.134 & (0.508, 0.988) & $\Delta\gamma$ & -0.153
& 0.155 & (-0.506, 0.106) \\
$\lambda_1$ & 0.784 & 0.130 & (0.531, 0.989) & & & & \\ \hline
%$\Delta$ & 0.019 & 0.186 & (-0.361, 0.393) \\
%$R$ & 1.067 & 0.326 & (0.508, 1.667) \\
%$\psi$ & 2.574 & 6.512 & (0.004, 8.772) \\
%$\delta$ & -0.135 & 0.176 & (-0.502, 0.221) \\
%$\Delta\gamma$ & -0.153 & 0.155 & (-0.506, 0.106) \\\hline
\multicolumn{8}{|c|}{{\bf Uniform Prior}}\\
\multicolumn{8}{|c|}{$P(\Delta>0|D)=0.583$, $P(\gamma_1-\gamma
_0<0|D)=0.854$, DIC$=$10.119, pD$=$1.352}\\ \hline
%& \multicolumn{3}{|c|}{$P(\gamma_1-\gamma_0<0|D)=0.854$}\\
%& \multicolumn{3}{|c|}{DIC=10.119,pD=1.352}\\\cline{2-4}
& mean & std & 95\% HPD & & mean & std & 95\% HPD \\\hline
$U$ & 0.833 & 0.141 & (0.548, 1.000) & $\Delta$ & 0.040 & 0.192 &
(-0.339, 0.418) \\
$V$ & 0.777 & 0.132 & (0.522, 0.989) & $R$ & 1.116 & 0.386 & (0.476,
1.831) \\
$\gamma_0$ & 0.231 & 0.148 & (0.007, 0.521) & $\psi$ & 2.161 & 3.727 &
(0.022, 6.814) \\
$\gamma_1$ & 0.075 & 0.073 & (0.000, 0.224) & $\delta$ & -0.110 & 0.186
& (-0.484, 0.267) \\
$\lambda_0$ & 0.686 & 0.140 & (0.414, 0.938) & $\Delta\gamma$ & -0.157
& 0.165 & (-0.515, 0.142) \\
$\lambda_1$ & 0.725 & 0.131 & (0.470, 0.956) & & & & \\ \hline
%$\Delta$ & 0.040 & 0.192 & (-0.339, 0.418)\\
%$R$ & 1.116 & 0.386 & (0.476, 1.831)\\
%$\psi$ & 2.161 & 3.727 & (0.022, 6.814)\\
%$\delta$ & -0.110 & 0.186 & (-0.484, 0.267) \\
%$\Delta\gamma$ & -0.157 & 0.165 & (-0.515, 0.142) \\\hline
\end{tabular}
\end{table}

\begin{table}[hbt]\caption{Number of scleroderma patients whose
forearm MRSS decreased
by 2 or 3, or has 0 MRSS at month 15.}\label{tab1sub2sec6}
\smallskip\centering%\footnotesize
\begin{tabular}{|lllcccc|}\hline
{Number of forearms} &&& \multicolumn{4}{c|}{Treatment Group}\\\cline{3-7}
{with improvement} &&& Collagen &&& Placebo\\\hline
0 &&& 36 &&&55\\
1 &&& 4 &&& 3\\
2 &&& 6 &&& 3\\
{\bf Total} &&& {\bf46} &&& {\bf61}\\\hline
\end{tabular}
\end{table}

\begin{table}[t!]\caption{Posterior estimates and credible intervals as
well as frequentist MLEs and confidence intervals
based on Dallal's model. }\label{tab2sub2sec6}
\smallskip\centering% \footnotesize
\begin{tabular}{|c|ccc|c|ccc|}\hline
\multicolumn{8}{|c|}{{\bf Jeffreys' Prior}} \\
\multicolumn{8}{|c|}{$P(\Delta>0|D)= 0.956$, DIC$=$18.563, pD$=$2.878,
$BF_\lambda=1.607$, $BF_\gamma=3.733$} \\ \hline
%& \multicolumn{3}{|c|}{$pD=2.878$} \\
%& \multicolumn{3}{|c|}{$BF_\lambda=1.607$, \; $BF_\gamma=3.733$} \\
& mean & std & 95\% HPD & & mean & std & 95\% HPD \\\hline
$U$ & 0.107 & 0.039 & (0.037, 0.186) & $\lambda_1$ & 0.178 & 0.049 &
(0.087, 0.277) \\
$V$ & 0.229 & 0.061 & (0.114, 0.348) & $\Delta$ & 0.095 & 0.057 &
(-0.017, 0.208) \\
$\gamma$ & 0.290 & 0.099 & (0.109, 0.486) & $R$ & 2.479 & 1.338 &
(0.640, 5.007) \\
$\lambda_0$ & 0.084 & 0.031 & (0.028, 0.147) & $\psi$ & 2.855 & 1.746 &
(0.566, 6.124) \\ \hline
%$\lambda_1$ & 0.178 & 0.049 & (0.087, 0.277) \\
%$\Delta$ & 0.095 & 0.057 & (-0.017, 0.208)\\
%$R$ & 2.479 & 1.338 & (0.640, 5.007) \\
%$\psi$ & 2.855 & 1.746 & (0.566, 6.124) \\\hline
\multicolumn{8}{|c|}{{\bf Bernardo's Prior}} \\
\multicolumn{8}{|c|}{$P(\Delta>0|D)=0.957$, DIC$=$18.523, pD$=$2.882,
$BF_\lambda=1.526$, $BF_\gamma=2.368$} \\ \hline
%& \multicolumn{3}{|c|}{$pD=2.882$} \\
%& \multicolumn{3}{|c|}{$BF_\lambda=1.526$, \; $BF_\gamma=2.368$}\\
& mean & std & 95\% HPD & & mean & std & 95\% HPD \\\hline
$U$ & 0.104 & 0.038 & (0.037, 0.182) & $\lambda_1$ & 0.174 & 0.049 &
(0.085, 0.272) \\
$V$ & 0.223 & 0.060 & (0.111, 0.341) & $\Delta$ & 0.092 & 0.056 &
(-0.015, 0.204) \\
$\gamma$ & 0.291 & 0.099 & (0.110, 0.487) & $R$ & 2.481 & 1.318 &
(0.630, 5.004) \\
$\lambda_0$ & 0.081 & 0.031 & (0.027, 0.143) & $\psi$ & 2.846 & 1.707 &
(0.582, 6.115) \\\hline
%$\lambda_1$ & 0.174 & 0.049 & (0.085, 0.272) \\
%$\Delta$ & 0.092 & 0.056 & (-0.015, 0.204)\\
%$R$ & 2.481 & 1.318 & (0.630, 5.004) \\
%$\psi$ & 2.846 & 1.707 & (0.582, 6.115) \\\hline
\multicolumn{8}{|c|}{{\bf Likelihood}}\\
\multicolumn{8}{|c|}{AIC=18.712, BIC$=$26.730, pD$=$3, $\chi^2_\lambda
=2.897$, $\chi^2_\gamma=0.152$}\\ \hline
%& \multicolumn{3}{|c|}{$pD=3$} \\
%& \multicolumn{3}{|c|}{$\chi^2_\lambda=2.897$, \; $\chi^2_
& mle & std & 95\% CI & & mle & std & 95\% CI \\\hline
$U$ & 0.098 & 0.038 & (0.024, 0.173) & $\lambda_1$ & 0.170 & 0.049 &
(0.073, 0.267) \\
$V$ & 0.217 & 0.061 & (0.098, 0.337) & $\Delta$ & 0.093 & 0.057 &
(-0.019, 0.205) \\
$\gamma$ & 0.280 & 0.102 & (0.081, 0.479) & $R$ & 2.210 & 1.057 &
(0.866, 5.641) \\
$\lambda_0$ & 0.077 & 0.030 & (0.017, 0.136) & $\psi$ & 2.458 & 1.323 &
(0.855, 7.062)\\\hline
%$\lambda_1$ & 0.170 & 0.049 & (0.073, 0.267)\\
%$\Delta$ & 0.093 & 0.057 & (-0.019, 0.205)\\
%$R$ & 2.210 & 1.057 & (0.866, 5.641)\\
%$\psi$ & 2.458 & 1.323 & (0.855, 7.062)\\\hline
\end{tabular}
\end{table}

Bayesian posterior estimates and 95\% HPD intervals for Dallal's model
and the saturated model
based on Jeffreys', the reference, and the uniform priors and 100,000
iterations are summarized
in Tables \ref{tab2sub1sec6} and \ref{tab3sub1sec6}. Under these three
priors, the posterior
means of $U$ are far away from the parameter space boundaries as
indicated by the non-zero
median estimates (not provided for obvious reasons). The 95\% HPD
intervals for $R$ under
the uniform prior tend to be larger than those under Jeffreys' prior or
the reference prior.
In addition, the uniform prior seems to perform the worst in terms of
the deviance information criterion (DIC) while Jeffreys' prior
and the reference prior seem to perform similarly. The risk difference
is essentially zero
($BF^{J}_\lambda=1.965$ and $P(\Delta>0|D)=0.550$; $BF^{R}_\lambda
=1.818$ and $P(\Delta>0|D)=0.541$) and the two cure rates
themselves are very high (above 75\%). We retain the saturated model
over Dallal's reduced model
according to the DIC and $P(\gamma_1-\gamma_0<0|D)$. However, the Bayes
factors indicate minimal evidence against the null hypothesis
$H^\ast_0: \gamma_1=\gamma_0$ ($BF^{J}_\gamma= 0.682$ and $BF^{R}_\gamma
=1.052$). The 95\% HPD interval for $\gamma_1-\gamma_0$
also covers zero. According to the retained model, the conditional
posterior probabilities of the cure rate at one site given the other
site was cured
are very high and slightly higher in the Amoxicillin group, therefore
this dependency cannot be ignored. There is also a noticeable discrepancy
between the correlation coefficients for the $Z_{ijk}$ variables for
the two treatment groups, $1- \dfrac{\gamma_i}{1-\lambda_i},\;i=0,1$.

\subsection{Bayesian Analysis of Bilateral Data with Large Sample
Data}\label{sub2sec6}

\citet{Postlethwaiteetal2008} %Postlethwaite et al. (2008)
considered a two-arm multi-centre double-blind randomized trial where
168 diffuse scleroderma patients are randomized to one of two groups to
receive either oral native collagen
at a dose of 500g/day or a similar appearing placebo. The total
duration of the treatment
phase was 12 months with an additional visit at month 15 for safety
follow-up. Rheumatologists routinely
examine both the left and right feet, forearms, hands, fingers, legs,
thighs, and upper arms of the
patient and assign a modified Rodnan Skin Score (MRSS) score between 0
and 3, that is 0 for normal, 1 for mild, 2 for moderate and 3
for severe skin thickening. The patient's improvement at each body part
level is recorded.
After consultation with rheumatologists, a patient has improved at a
body part level if the MRSS at the
end of the trial is either zero or has dropped by two units or more
from baseline. The goal is to test
whether there is a significant difference in the improvement rates
between the two groups at each body
part level. Table \ref{tab1sub2sec6} reports an examplary set of
results from the trial for the forearms.

In blocks 1 and 2 of Table \ref{tab2sub2sec6}, we provide
Bayesian posterior estimates and credible intervals for Dallal's model
based on Jeffreys' and Bernardo's priors. The results are based on
100,000 posterior simulations.
In block 3 of Table \ref{tab2sub2sec6}, we provide frequentist MLEs,
Wald FCIs, Aikaike information criterion (AIC) and
Bayesian information criterion (BIC), along with large-sample $\chi^2$
test for $H_0:\lambda_0=\lambda_1$
and $H^\ast_0: \gamma_1=\gamma_0$. For this data, the total sample size in
each group is moderate. Although four of the cells have counts less
than 5, the tallies that matter here are the
subtotals, $m_{1+}=7, m_{2+}=9, m_{10}+m_{20}=6, m_{00}=55,
m_{11}+m_{21}=10, m_{01}=36$, which
are all greater than 5. We provide only the results for Dallal's
reduced model.
Indeed, Bayes factors indicate the null hypothesis, $H^\ast_0: \gamma
_1=\gamma_0$, is supported.
In this example, the treatment cure rate has decreased the disease risk
by two fold. The risk of disease in both
treatment and control groups remain high (above 80\%). The null
hypothesis of equality of the treatment cure rates,
$H_0:\lambda_0=\lambda_1$ is supported by our analysis. Overall
Dallal's model and the saturated
model give similar results and there is no difference between using
Jeffreys' prior (e.g., $\widehat{R}=2.479$ and BCI:
(0.640, 5.007)) or Bernardo's prior (e.g., $\widehat{R}=2.481$ and
BCI: (0.630, 5.004)), although results from Bernardo's
reference prior are easier to compute. Note that our posterior
estimates and HPD intervals for the risk difference are
in line with the results in Pei et al. (2010) obtained under the equal
correlation model
($\widehat{\Delta}=0.0970$, CI: $(-0.0214, 0.2217)$).

%%----------------------------------------------------------------------------------------------
% Section 8
%%----------------------------------------------------------------------------------------------

\section{General Classes of Prior Distributions}\label{sec8}

The reference priors and the uniform prior discussed earlier can be
embedded in the family of prior distributions\vadjust{\eject}
\begin{eqnarray}\label{eq4sub3sec4}
\pi(\gamma, \lambda_0, \lambda_1) & = &
\dfrac{2^\alpha}{\B(\alpha,\beta)} \dfrac{\gamma^{\alpha-1}
(1-\gamma)^{\beta-1}}{(1+\gamma)^{\alpha+\beta-a_0-a_1}}
\dfrac{\lambda_0^{a_0-1} \big[1-(1+\gamma)\lambda_0\big]^{b_0-1}}{\B
(a_0,b_0)}\nonumber\\
& & \dfrac{\lambda_1^{a_1-1}
\;\big[1-(1+\gamma)\lambda_1\big]^{b_1-1}}{\B(a_1,b_1)},\qquad
(\gamma, \lambda_0, \lambda_1) \ \in\Omega,
\end{eqnarray}
which is equivalent to stating that
\begin{equation}\label{eq5sub3sec4}
\pi(\gamma, u, v) = \dfrac{2^\alpha}{\B(\alpha,\beta)}
\dfrac{\gamma^{\alpha-1} (1-\gamma)^{\beta-1}}{(1+\gamma)^{\alpha+\beta}}
\dfrac{u^{a_0-1}\;(1-u)^{b_0-1}}{\B(a_0,b_0)} \dfrac{v^{a_1-1}
(1-v)^{b_1-1}}{\B(a_1,b_1)}.
\end{equation}
Another representation of this class of prior distributions is through
the following hierarchical model:
(a) $\gamma\sim f(\gamma) = \dfrac{2^\alpha}{\B(\alpha,\beta)}\;
\dfrac{\gamma^{\alpha-1}\;(1-\gamma)^{\beta-1}}{(1+\gamma)^{\alpha+\beta
}}$ or equivalently
$\dfrac{1-\gamma}{1+\gamma}\sim\Be(\beta, \alpha)$,
(b) $\lambda_0|\gamma\sim\Be\left(a_0, b_0; 0,\dfrac{1}{1+\gamma
}\right)$,
and (c) $\lambda_1|\gamma\sim\Be\left(a_1, b_1; 0,\dfrac{1}{1+\gamma
}\right)$,
where $\Be(\alpha,\beta; l, u)$ stands for a Beta random variable with
shape parameters $\alpha$
and $\beta$ defined on the interval $(l,u)$.

Jeffreys' prior distribution for the parameterization $(\gamma, U, V)$
suggests a
larger family of conjugate prior distributions for $(\gamma, U, V)$, namely,
\[
\pi(\gamma, u, v) \propto
\frac{2^\alpha}{\B(\alpha,\beta)}\frac{\gamma^{\alpha-1}
(1-\gamma)^{\beta-1}}{(1+\gamma)^{\alpha+\beta}} (u + rv)^{1/2}
u^{a_0-1}(1-u)^{b_0-1} v^{a_1-1}(1-v)^{b_1-1},
\]
with $0<\gamma, u, v<1$ and $\alpha,\beta, a_0, b_0, a_1, b_1>0$, which
translates into the prior distribution%\vadjust{\eject}
\begin{align}\label{eq7sub2sec4}
\pi(\gamma, \lambda_0, \lambda_1) \propto&
\; \frac{2^\alpha}{\B(\alpha,\beta)}\;\frac{\gamma^{\alpha-1}
(1-\gamma)^{\beta-1}}{(1+\gamma)^{\alpha+\beta-a_0-a_1-d}} (\lambda_0 +
r\lambda_1)^{1/2}
\lambda_0^{a_0-1}\big[1-(1+\gamma)\lambda_0\big]^{b_0-1}\nonumber\\
& \; \times\lambda_1^{a_1-1}\big[1-(1+\gamma)\lambda_1\big
]^{b_1-1},\qquad\qquad(\gamma,\lambda_0,\lambda_1)\in\Omega.
\end{align}

\section{Concluding Remarks}\label{sec7}

Using the parameterization $(\gamma,U,V)$, it can be deduced that
$\Delta=\dfrac{V-U}{1+\gamma}$.
This result highlights a direct dependence of the risk difference on
the nuisance parameter $\gamma$.
This result also points out the main difference between the risk
difference in $3\times2$ bilateral data
and the risk difference, $V-U$, in an ordinary $2\times2$ table where
one collects a single measurement
per subject. In other words, the divisor $1+\gamma$ is the term that
connects the two disease risks in
the bilateral data context. Indeed, according to the expression of the
likelihood, $L(\gamma, U, V)$,
given in Appendix 1, %\ref{sec2app},
the parameters $U$ and $V$ can be interpreted as the proportions of
cases with one or more body part(s) cured in the placebo and treatment
groups. On the opposite side,
$\lambda_0$ and $\lambda_1$ are interpreted as the proportions of body
parts cured in the placebo and
treatment groups. So unlike $U$ and $V$, one individual can contribute
twice in the computation of
$\lambda_0$ and $\lambda_1$.

\begin{table}[hbpt]\caption{Relevant summary statistics for the risk
ratio.}\label{tab1sec7}
\centering
\begin{tabular}{lllll}\hline
& \multicolumn{4}{c}{Group} \\\cline{2-5}
Numbers of & & & & \\
cured organs & Treatment & & & Control\\\hline
0 & $m_{01} \; (1-U)$ & & & $m_{00} \; (1-V)$\\
1 or 2 & $m_{+1}-m_{01} \; (U)$ & & & $m_{+0}-m_{00} \; (V)$\\\hline
\end{tabular}
\end{table}

Another benefit of the parameterization $(\gamma,U,V)$ is that it shows
that the risk ratio,
$R=\dfrac{\lambda_1}{\lambda_0}=\dfrac{V}{U}$, does not depend on the
nuisance parameter
$\gamma$. As a result, both frequentist and Bayesian inferences do not
depend on $\gamma$ and are easier to compute.
Therefore,
the risk ratio has a technical advantage on the risk difference.
Moreover, the definition of $R$ in
$3\times2$ bilateral data coincides with the definition $R$ from the
$2\times2$ binary table in Table \ref{tab1sec7}.
Actually, the $3\times2$ bilateral table can be replaced by Table \ref
{tab1sec7} when the focus is on $R$.
For these reasons, our choice of the parameter in a bilateral data
design is the risk ratio.

%%----------------------------------------------------------------------------------------------
% Section 7
%%----------------------------------------------------------------------------------------------

While there are numerous frequentist papers dealing with bilateral
data, this work remains incomplete.
For example, Dallal's model used in this paper has not been
investigated in the frequentist
literature although other models have and there are no Bayesian
treatments of the problem.
Although the risk ratio and the odds ratio are well established
parameters in medical settings,
they do not appear in the bilateral data literature. A clear advantage
of the risk ratio over the commonly used
risk difference is that inference does not involve the nuisance
parameter $\gamma$.
In addition, the risk ratio remains unchanged when going from an
ordinary $2\times2$ binary table where
all observations are independent to a $3\times2$ binary table where
observations taken from the same subjects are
correlated. The risk ratio and the odds ratio open the door for
bilateral data to be studied under a prospective scheme
or a retrospective scheme using an observational study paradigm.
The presentation exposed here has taken into account all these
inconveniences and has provided a broad
discussion of the Bayesian framework both from the point of view of the
tests of hypotheses as well as the estimation
of the key model parameters. We have added a simulation study to
empirically compare the effectiveness
of Bayesian methods against themselves as well as frequentist methods
in the context of small, moderate
and large sample sizes and for a wide range of $\Delta$ values. For
example, we have found that Jeffreys'
prior and the reference prior tend to perform similarly. We have also
found that frequentist methods tend to perform
very poorly when $\Delta$ is small and the sample size is small or
moderate. The uniform prior has the best
overall property when it comes to estimating the parameter $\gamma$. We
have concluded our work with two detailed case studies,
one of which shows that it is impossible to carry out frequentist
inference given that some of the parameters sit on
the boundary of the parameter space or the normal approximation is not accurate.
Our Bayesian framework works remarkably well in these situations as
well as the large sample cases.

When subject level bilateral data with covariates are available,
Dallal's regression model can be developed in order to incorporate covariates.
As discussed in Section \ref{sub2sec4}, $\gamma$, $U = (1+\gamma)\lambda
_0$, and $V = (1+\gamma)\lambda_1$ are unconstrained
and the parameter space is (0,1) for each of these three transformed
parameters. Therefore,
a logistic regression model can be assumed for each of $\gamma$, $U$,
and $V$. Consequently, Jeffreys' prior and the reference prior can be derived.
However, the computational and theoretical properties of these priors
need to be carefully examined. The development of
Dallal's regression model deserves a future research project, which is
currently under investigation.

%
% Appendix
%

\begin{acknowledgement}

The authors wish to thank the Editor-in-Chief, the Editor, the Associate
Editor, and the anonymous referee for their very helpful comments and
suggestions,
which have led to a much improved version of the paper. Dr. M.-H.
Chen's research was partially
supported by NIH grants \#GM 70335 and \#CA 74015.%\vadjust{\vfill{\eject}}

\end{acknowledgement}

%% Supplement Material %%
\begin{supplement}
% \sname{}\label{}
\stitle{Supplementary Web Materials for ``Objective Bayesian Inference for Bilateral Data''}
\slink[doi]{10.1214/14-BA890SUPP}
\sdatatype{.pdf}
\sfilename{mainAppendixBilateralStudiesPaperR2-Jun-27-2014.pdf}
\end{supplement}

\bibliographystyle{ba}

\begin{thebibliography}{24}
\newcommand{\enquote}[1]{``#1''}
\expandafter\ifx\csname natexlab\endcsname\relax\def\natexlab#1{#1}\fi
\expandafter\ifx\csname url\endcsname\relax
\def\url#1{\texttt{#1}}\fi
\expandafter\ifx\csname urlprefix\endcsname\relax\def\urlprefix{URL }\fi

\bibitem[{Berger and Bernardo(1989)}]{BergerBernardo1989}
Berger, J.~O. and Bernardo, J.~M. (1989).
\newblock\enquote{Estimating a Product of Means: {Bayesian} Analysis with
Reference Priors.}
\newblock\emph{Journal of the American of Statistical Assocation\/}, 84:
200--207.
\endbibitem

\bibitem[{Berger and Bernardo(1992a)}]{BergerBernardo1992a}
--- (1992a).
\newblock\enquote{Ordered Group Reference Priors with Applications to a
Multinomial Problem.}
\newblock\emph{Biometrika\/}, 79: 25--37.
\endbibitem

\bibitem[{Berger and Bernardo(1992b)}]{BergerBernardo1992b}
--- (1992b).
\newblock\enquote{Reference Priors in a Variance Components Problem.}
\newblock In Goel, P.~K. and Iyengar, N. (eds.), \emph{Bayesian Analysis in
Statistics and Econometrics\/}, 323--340. New York: Springer-Verlag.
\endbibitem

\bibitem[{Berger and Bernardo(1992c)}]{BergerBernardo1992c}
--- (1992c).
\newblock\enquote{On the Development of Reference Priors.}
\newblock In Bernado, J.~M., Berger, J.~O., Dawid, A.~P., and Smith, A. F.~M.
(eds.), \emph{Bayesian Statistics\/}, volume~4, 35--60. New York: Oxford:
University Press.
\endbibitem

\bibitem[{Bernardo(1981)}]{Bernardo1979}
Bernardo, J.~M. (1981).
\newblock\enquote{Reference Posterior Distributions for {Bayes} Inference.}
\newblock\emph{Journal of the Royal Statistical Society, Series B\/}, 41:
113--147.
\endbibitem

\bibitem[{Chen and Shao(1999)}]{ChenShao1999}
Chen, M.-H. and Shao, Q.-M. (1999).
\newblock\enquote{{Monte Carlo} Estimation of {Bayesian} Credible and {HPD}
Intervals.}
\newblock\emph{Journal of Computational and Graphical Statistics\/}, 8: 69--92.
\endbibitem

\bibitem[{Cox and Reid(1987)}]{CoxReid1987}
Cox, D.~R. and Reid, N. (1987).
\newblock\enquote{Parameter Orthogonality and Approximate Conditional
Inference.}
\newblock\emph{Journal of the Royal Statistical Society, Series B\/}, 49:
1--39.
\endbibitem

\bibitem[{Dallal(1988)}]{Dallal1988}
Dallal, G.~E. (1988).
\newblock\enquote{Paired Bernoulli Trials.}
\newblock\emph{Biometrics\/}, 44: 253--257.
\endbibitem

\bibitem[{Datta and Ghosh(1996)}]{DattaGhosh1996}
Datta, G.~S. and Ghosh, M. (1996).
\newblock\enquote{On the Invariance of Noninformative Priors.}
\newblock\emph{The Annals of Statistics\/}, 24: 141--159.
\endbibitem

\bibitem[{Ghosh and Mukerjee(1992)}]{GhoshMukerjee1992}
Ghosh, J.~K. and Mukerjee, R. (1992).
\newblock\enquote{Non-informative Priors.}
\newblock In Bernardo, J.~M., Berger, J.~O., Dawid, A.~P., and Smith,
A. F.~M.
(eds.), \emph{Bayesian Statistics\/}, volume~4, 195--210. Oxford: Oxford
University Press.
\endbibitem

\bibitem[{Jeffreys(1946)}]{Jeffreys1946}
Jeffreys, H. (1946).
\newblock\enquote{An Invariant Form for the Prior Probability in Estimation
Problems.}
\newblock\emph{Proceedings of the Royal Society of London. Series A,
Mathematical and Physical Sciences\/}, 186: 453--461.
\endbibitem

\bibitem[{Mandel et~al.(1982)Mandel, Bluestone, Rockette, Blatter, Reisinger,
Wucher, and Harper}]{MandelBluestone1982}
Mandel, E.~M., Bluestone, C.~D., Rockette, H.~E., Blatter, M.~M., Reisinger,
K.~S., Wucher, F.~P., and Harper, J. (1982).
\newblock\enquote{Duration of Effusion after Antibiotic Treatment for Acute
Otitis Media: Comparison of Cefaclor and Amoxicillin.}
\newblock\emph{Pediatric Infectious Diseasee\/}, 1: 310--316.
\endbibitem

\bibitem[{Morris(1993)}]{Morris1993}
Morris, R.~W. (1993).
\newblock\enquote{Bilateral Procedures in Randomised Controlled Trials.}
\newblock\emph{The Journal of Bone and Joint Surgery\/}, 75: 675--6.
\endbibitem

\bibitem[{Pei et~al.(2008)Pei, Tang, and Guo}]{PeiTangGuo2008}
Pei, Y.-B., Tang, M.-L., and Guo, J. (2008).
\newblock\enquote{Testing the Equality of Two Proportions for Combined
Unilateral and Bilateral Data.}
\newblock\emph{Communications in Statistics -- Simulation and
Computation\/},
37: 1515--1529.
\endbibitem

\bibitem[{Pei et~al.(2010)Pei, Tang, Wong, and Guo}]{PeiTangWongGuo2010}
Pei, Y.-B., Tang, M.-L., Wong, W.-K., and Guo, J. (2010).
\newblock\enquote{Confidence Intervals for Correlated Proportion Differences
from Paired Data in a Two-Arm Randomized Clinical Trial.}
\newblock\emph{Statistical Methods in Medical Research\/}, 21: 167--187.
\endbibitem

\bibitem[{Postlethwaite et~al.(2008)Postlethwaite, Wong, Clements, Chatterjee,
Fessler, Kang, Korn, Mayes, Merkel, Molitor, Moreland, Rothfield, Simms,
Smith, Spiera, Steen, Warrington, White, Wigley, and
Furst}]{Postlethwaiteetal2008}
Postlethwaite, A.~E., Wong, W.~K., Clements, P., Chatterjee, S., Fessler,
B.~J., Kang, A.~H., Korn, J., Mayes, M., Merkel, P.~A., Molitor, J.~A.,
Moreland, L., Rothfield, N., Simms, R.~W., Smith, E.~A., Spiera, R., Steen,
V., Warrington, K., White, B., Wigley, F., and Furst, D.~E. (2008).
\newblock\enquote{A Multicenter, Randomised, Double-Blind, Placebo-Controlled
Trial of Oral Type I Collagen in Patients with Diffuse Cutaneous Systemic
Sclerosis: I. Oral Type I Collagen Does not Improve Skin in all
Patients, but
may Improve Skin in Late-Phase Disease.}
\newblock\emph{Arthritis and Rheumatism\/}, 58: 1810--1822.
\endbibitem

\bibitem[{Qiu et~al.(2009)Qiu, Tang, and Tang}]{QiuTangTang2009}
Qiu, S.-F., Tang, N.-S., and Tang, M.-L. (2009).
\newblock\enquote{Sample Size for Testing Difference between two Proportions
for the Bilateral-Sample Design.}
\newblock\emph{Journal of Biopharmaceutical Statistics\/}, 19: 857--871.
\endbibitem

\bibitem[{Rosner(1982)}]{Rosner1982}
Rosner, B. (1982).
\newblock\enquote{Statistical Methods in Ophthalmology: An Adjustment
for the
Intraclass Correlation between Eyes.}
\newblock\emph{Biometrics\/}, 38: 105--114.
\endbibitem

\bibitem[{Sun and Berger(1998)}]{SunBerger1998}
Sun, D. and Berger, J.~O. (1998).
\newblock\enquote{Reference Priors with Partial Information.}
\newblock\emph{Biometrika\/}, 85: 55--71.
\endbibitem

\bibitem[{Tang et~al.(2010)Tang, Pei, Wong, and Li}]{TangPeiWongLi2010}
Tang, M.-L., Pei, Y.-B., Wong, W.-K., and Li, J.-L. (2010).
\newblock\enquote{Goodness-of-fit Tests for Correlated Paired Binary Data.}
\newblock\emph{Statistical Methods in Medical Research\/}, 1--15.
\endbibitem

\bibitem[{Tang et~al.(2006)Tang, Tang, and Rosner}]{TangTangRosner2006}
Tang, M.-L., Tang, N.-S., and Rosner, B. (2006).
\newblock\enquote{Statistical Inference for Correlated Data in Ophthalmologic
Studies.}
\newblock\emph{Statistics in Medicine\/}, 25: 2771--2783.
\endbibitem

\bibitem[{Tang et~al.(2011)Tang, Qui, Tang, and Pei}]{TangQiuTangPei2011}
Tang, N.-S., Qui, S.-F., Tang, M.-L., and Pei, Y.-B. (2011).
\newblock\enquote{Asymptotic Confidence Interval Construction for Proportion
Difference in Medical Studies with Bilateral Data.}
\newblock\emph{Statistical Methods in Medical Research\/}, 20:
233--259.%\vadjust{\eject}
\endbibitem

\bibitem[{Tang et~al.(2008)Tang, Tang, and Qiu}]{TangTangQiu2008}
Tang, N.-S., Tang, M.-L., and Qiu, S.-F. (2008).
\newblock\enquote{Testing the Equality of Proportions for Correlated
Otolaryngologic Data.}
\newblock\emph{Computational Statistics and Data Analysis\/}, 52: 3719--3729.
\endbibitem

\bibitem[{Yang(1995)}]{Yang1995}
Yang, R. (1995).
\newblock\enquote{Invariance of the Reference Prior Under Reparametrization.}
\newblock\emph{Test\/}, 4: 83--94.
\endbibitem

\end{thebibliography}

%% Appendix %%
\appendix
% \section{}\label{}

\setcounter{equation}{0}
\def\theequation{A.\arabic{equation}}

\section*{Appendix 1: Derivation of Jeffreys' Prior}\label{sec2app}

Set $U = (1+\gamma)\lambda_0$ and $V = (1+\gamma)\lambda_1$. Under this new
parametrization, the parameter space reduces to the interval $[0,1]$ for
each of the three parameters. Theorem A.1 %\ref{th1sec2app}
emphasizes
the role of a reparametrization technique in obtaining Jeffreys' prior
distribution
in the original parametrization after having derived Jeffreys' prior in
the reparametrized
space.

\medskip

\noindent
{\bf Theorem A.1} (Consistency Under Reparametrization). {\it
Consider a model \linebreak
$\mathfrak{M}\equiv\left\{p(X|\theta), \x\in\X,\theta\in\Theta\right\}
$ and let
$\phi(\theta)$ be an invertible transformation of $\theta$. Then, the
Jeffreys' prior
corresponding to the parameter $\phi$, $\pi(\phi)$, is that induced by
the Jeffreys'
prior density of $\theta$, $\pi(\theta)$.
}
We prove Proposition \ref{prop1sub2sec4} and \ref{prop2sub2sec4} below.
\begin{proof}
As a function of $(\gamma, U,V)$, the likelihood function simplifies to
\[
L(\gamma, U, V) \propto\dfrac{\gamma^{m_{1+}}
(1-\gamma)^{m_{2+}}}{(1+\gamma)^{m_{1+}+m_{2+}}}
U^{m_{10}+m_{20}}(1-U)^{m_{00}} V^{m_{11}+m_{21}}
(1-V)^{m_{01}}.
\]
This parametrization splits the likelihood function into three
unrelated pieces, each piece
related to a single parameter. This makes it very easy to derive
Jeffreys' prior.
The first and second derivatives of the log-likelihood function with
respect to $\gamma$, $U$, and $V$ are
\begin{align*}
\frac{\partial l(\gamma)}{\partial\gamma} = \dfrac{m_{1+}}{\gamma}-
\dfrac{m_{2+}}{1-\gamma} -
\dfrac{m_{1+}+m_{2+}}{1+\gamma}, & \; \frac{\partial^2 l(\gamma
)}{\partial\gamma^2} = -
\dfrac{m_{1+}}{\gamma^2}-\dfrac{m_{2+}}{(1-\gamma)^2}+\dfrac
{m_{1+}+m_{2+}}{(1+\gamma)^2}, \\
\frac{\partial l(U)}{\partial U} = \dfrac{(m_{20} +
m_{10})}{U}-\dfrac{m_{00}}{1-U}, & \; \frac{\partial^2 l(U)}{\partial
U^2} = -\dfrac{(m_{20} +
m_{10})}{U^2} - \dfrac{m_{00}}{(1-U)^2}, \\
\frac{\partial l(V)}{\partial V} = \dfrac{(m_{21}+m_{11})}{V}\;-\;\dfrac
{m_{01}}{1-V}, & \;
\frac{\partial^2 l(V)}{\partial V^2} = -\dfrac
{(m_{21}+m_{11})}{V^2}-\dfrac{m_{01}}{(1-V)^2}.
\end{align*}
Let $r=\dfrac{m_{+1}}{m_{+0}}$ be the ratio of the sample size in the
two groups. Thus, we have
\begin{align*}
E\Big[-\frac{\partial^2 l(\gamma)}{\partial\gamma^2}\Big] = & \;
\dfrac{2m_{+0}(U + rV)}{\gamma(1-\gamma)(1+\gamma)^2}, \;\;
E\Big[-\dfrac{\partial^2 l(U)}{\partial U^2}\Big] =
\dfrac{m_{+0}}{U(1-U)}, \\
E\Big[-\frac{\partial^2 l(V)}{\partial V^2}\Big] = & \;
\dfrac{m_{+1}}{V(1-V)}.
\end{align*}
Hence, Jeffreys' prior under the parameterization $(\gamma, U, V)$ is
\[
\pi_J(\gamma, u, v) \propto \sqrt{\dfrac{(u + rv) }{\gamma
(1-\gamma) (1+\gamma)^2 u(1-u) v(1-v)}}, \qquad0<\gamma, u, v<1,
\]
and it is proper. Theorem A.1 %\ref{th1sec2app}
implies that
Jeffreys' prior density in the parameterization $(\gamma, \lambda_0,
\lambda_1)$ is
\[
\pi_J(\gamma, \lambda_0, \lambda_1) \propto
\sqrt{\dfrac{(1+\gamma)(\lambda_0 + r\lambda_1) }{\gamma
(1-\gamma) \lambda_0 \lambda_1
\big[1-(1+\gamma)\lambda_0\big]\; \big[1-(1+\gamma)\lambda_1\big]}},
\]
where $(\gamma, \lambda_0, \lambda_1) \in\Omega$.%\vadjust{\eject}
\end{proof}

\section*{Appendix 2: Derivation of Reference Priors}\label{sec1app}

\begin{proof}
We first derive the joint reference prior for the parameterization
$(\gamma, U, V)$ and
then transform back to derive the joint reference prior for $(\gamma,
\lambda_0, \lambda_1)$.
See \citet{Yang1995} %Yang (1995)
for justification. We also adopt the notation in \citet{Yang1995}. %
%Yang (1995).

\nd{\bf Case 1:} Consider the group ordering $\{U,V\}$ and then $\gamma$
or $\{V,U\}$ and then $\gamma$. We have
$h_1=\dfrac{m_{+0}}{u(1-u)}$, $h_2=\dfrac{m_{+1}}{v(1-v)}$, and
$h_3=\dfrac{2m_{+0}(u+rv)}{\gamma(1-\gamma)(1+\gamma)^2}$. Thus, we obtain
\begin{align*}
\pi^1_R(u,v,\gamma) = & \; \dfrac{|h_3|^{1/2}}{\int|h_3|^{1/2}
d\gamma}\dfrac{\exp\left\{\frac{1}{2}
\int\log(h_1h_2)\,d\gamma\right\}} {\exp\left\{\dfrac{1}{2}\, \iiint
\log(h_1h_2) d\gamma du dv\right\}},\\
\propto& \; \dfrac{\sqrt{2}}{\pi} \dfrac{\gamma^{1/2-1}
(1-\gamma)^{1/2-1}}{(1+\gamma)}
\dfrac{u^{1/2-1}(1-u)^{1/2-1}}{\B(1/2,1/2)}\dfrac
{v^{1/2-1}(1-v)^{1/2-1}}{\B(1/2,1/2)},
\end{align*}
and
\begin{align*}
& \; \pi^1_R(\gamma, \lambda_0, \lambda_1) \\
= & \; \dfrac{\sqrt{2}}{\pi} \gamma^{1/2-1}
(1-\gamma)^{1/2-1}\; \dfrac{\lambda_0^{1/2-1}
\big[1-(1+\gamma)\lambda_0\big]^{1/2-1}}{\B(1/2,1/2)}
\dfrac{\lambda_1^{1/2-1}\big[1-(1+\gamma)\lambda_1\big]^{1/2-1}}{\B(1/2,1/2)}.
\end{align*}

\nd{\bf Case 2:} Consider the group ordering of $\{\gamma\}$ and then
$\{U,V\}$ or the ordering $\{\gamma\}$
and then $\{V,U\}$. We have
$\,h_1=\dfrac{2m_{+0}(u+rv)}{\gamma(1-\gamma)(1+\gamma)^2},\,
h_2=\dfrac{m_{+0}}{u(1-u)}, \, h_3=\dfrac{m_{+1}}{v(1-v)}$.
Thus, the reverse reference prior is
\begin{align*}
\pi^2_R(u,v,\gamma) = & \; \dfrac{|h_1\,h_2|^{1/2}}{\int|h_1 h_2|^{1/2}
dudv}\dfrac{\exp\left\{\frac{1}{2}
\iint\log(h_3)\,du\,dv\right\}} {\exp\left\{\frac{1}{2} \iiint
\log(h_3) du dv d\gamma\right\}}\nonumber\\
= & \; \dfrac{u^{1/2-1}(1-u)^{1/2-1}}{\B(1/2,1/2)} \dfrac
{v^{1/2-1}(1-v)^{1/2-1}}{\B(1/2,1/2)}
\dfrac{\sqrt{2}}{\pi} \dfrac{\gamma^{1/2-1} (1-\gamma
)^{1/2-1}}{(1+\gamma)} \\
= & \pi^1_R(u,v,\gamma).
\end{align*}\vadjust{\eject}

\nd{\bf Case 3:} Here one starts with some subjective joint prior
distribution for the parameters for which one has
a good knowledge of and for the other parameters, one uses a
non-informative prior
distribution to reflect the lack of knowledge. Assume there is prior
evidence for assuming the following conditional joint distribution
\[
\pi^3_B(\lambda_0, \lambda_1| \gamma) = \dfrac{\lambda_0^{a_0-1}
\left(\frac{1}{1+\gamma}-\lambda_0\right)^{b_0-1}}
{\left(\frac{1}{1+\gamma}\right)^{a_0+b_0-1}\B(a_0,b_0)}
\dfrac{\lambda_1^{a_1-1}
\left(\frac{1}{1+\gamma}-\lambda_1\right)^{b_1-1}}{\left(\frac
{1}{1+\gamma}\right)^{a_1+b_1-1}\B(a_1,b_1)},\quad
0<\lambda_0, \lambda_1<\frac{1}{1+\gamma}.
\]
about $\lambda_0$ and $\lambda_1$ given $\gamma$. That is,
\[
\lambda_0|\gamma\sim\Be\left(a_0, b_0; 0,\dfrac{1}{1+\gamma}\right),
\;
\lambda_1|\gamma\sim\Be\left(a_1,b_1; 0, \dfrac{1}{1+\gamma}\right)
\]
and both $\lambda_0$ and
$\lambda_1$ are conditionally independent given $\gamma$. But one has
no idea about a prior for $\gamma$.
A solution to this problem is to use a reference prior for $\gamma$.
Under our partial prior specification, the
reference prior is Jeffreys' prior associated to the integrated
likelihood (integrating out $U$ and $V$),
\[
L(\gamma) \propto \dfrac{\gamma^{m_{1+}}
(1-\gamma)^{m_{2+}}}{(1+\gamma)^{m_{1+}+m_{2+}}}.
\]
Hence, the reference prior corresponding to this integrated
likelihood is
\[
\pi^3_R(\gamma) \propto \sqrt{\dfrac{m_{+0}E(\lambda_0) +
m_{+1}E(\lambda_1) }{\gamma(1-\gamma) (1+\gamma)^2}} \propto
\sqrt{\dfrac{1}{\gamma(1-\gamma) (1+\gamma)^2}},
\]
and it is proper. Thus, the joint prior distribution over $\Omega$ is then
\begin{align*}
& \; \pi^3_R(\gamma, \lambda_0, \lambda_1) \\
= & \;
\dfrac{\sqrt{2}}{\pi}\dfrac{\gamma^{1/2-1}
(1-\gamma)^{1/2-1}}{(1+\gamma)} \dfrac{\lambda_0^{a_0-1}
\;\left(\frac{1}{1+\gamma}-\lambda_0\right)^{b_0-1}}
{\left(\frac{1}{1+\gamma}\right)^{a_0+b_0-1}\B(a_0,b_0)}\dfrac{\lambda_1^{a_1-1}
\left(\frac{1}{1+\gamma}-\lambda_1\right)^{b_1-1}}{\left(\frac
{1}{1+\gamma}\right)^{a_1+b_1-1}\B(a_1,b_1)},
\end{align*}
and it belongs to the family of prior distributions discussed in \eqref
{eq4sub3sec4}.

\nd{\bf Case 4:} Similarly, in the second setup of partial prior
specification, one
assumes that there is available a family of prior distributions for
$\gamma$, for example,
\[
\pi^4_B(\gamma) = \dfrac{2^\mu}{\B(\mu,\nu)}\dfrac{\gamma^{\mu-1}
(1-\gamma)^{\nu-1}}{(1+\gamma)^{\mu+\nu}},
\]
and one would like to find the joint reference prior distribution for
the pair
$(\lambda_0, \lambda_1)$ conditional on $\gamma$. The proposed
reference prior
under this partial prior specification is
\[
\pi^4_R(\lambda_0, \lambda_1|\gamma) \propto \big|\det(S)\big|^{1/2} =
\sqrt{F_{22} F_{33}} \propto\dfrac{(1+\gamma)}{\sqrt{\lambda_0\lambda_1
\big[1-(1+\gamma)\lambda_0\big]\big[1-(1+\gamma)\lambda_1\big]}},
\]
where $S$ is the lower $2\times2$ left corner matrix of the Fisher
information matrix. Hence, the joint prior distribution of
$(\gamma,\lambda_0,\lambda_1)$ over $\Omega$ is
\begin{align*}
& \; \pi^4_R(\gamma, \lambda_0, \lambda_1) \\
= & \;
\dfrac{2^\mu}{\B(\mu,\nu)}\dfrac{\gamma^{\mu-1}
(1-\gamma)^{\nu-1}}{(1+\gamma)^{\mu+\nu-1}}\dfrac{\lambda_0^{1/2-1}
\left[1-(1+\gamma)\lambda_0\right]^{1/2-1}}{\B(1/2,1/2)}\dfrac{\lambda_1^{1/2-1}
\left[1-(1+\gamma)\lambda_1\right]^{1/2-1}}{\B(1/2,1/2)}
\end{align*}
and it belongs again to the family of prior distributions discussed in
\eqref{eq4sub3sec4}.
\end{proof}

\section*{Appendix 3: Results for Posterior Calculation}\label{sec3app}

\noindent
{\bf Proposition A.1.} {\it
Let $\gamma$ have density $f(\gamma) =
\dfrac{2^{\mu}}{\B(\mu,\nu)}\dfrac{\gamma^{\mu-1} (1-\gamma)^{\nu-1}}{
(1+\gamma)^{\mu+\nu}}, \; 0<\gamma<1$. Then, the density of $\phi=\logit
(\gamma)+\log(2)$
is $ f(\phi) = \dfrac{1}{\B(\mu,\nu)}\dfrac{e^{\mu\phi}}{\big(1+e^\phi
\big)^{\mu+\nu}}$,
which is well known to be the density of $\logit(p)$, where $p\sim\Be
(\mu,\nu)$.}
\begin{proof}
$\gamma= \dfrac{e^{\phi-\log(2)}}{1+e^{\phi-\log(2)}}$.
We have $1-\gamma= \dfrac{1}{1+e^{\phi-\log(2)}}$, $1+\gamma=
\dfrac{1+2e^{\phi-\log(2)}}{1+e^{\phi-\log(2)}}$, and $\dfrac{d\gamma
}{d\phi} =
\dfrac{e^{\phi-\log(2)}}{(1+e^{\phi-\log(2)})^2}$. Thus, we obtain
\begin{align*}
\hspace*{-20pt}f(\phi) = & \;
\dfrac{2^\mu}{\B(\mu,\nu)}
\dfrac{e^{(\mu-1){(\phi-\log(2))}}}{(1+e^{\phi-\log(2)})^{\mu-1}} \dfrac
{1}{(1+e^{\phi-\log(2)})^{\nu-1}}
\dfrac{(1+e^{\phi-\log(2)})^{\mu+\nu}}{
(1+2e^{\phi-\log(2)})^{\mu+\nu}} \dfrac{e^{\phi-\log(2)}}{(1+e^{\phi
-\log(2)})^2}\hspace*{-20pt}\\
= & \; \dfrac{2^\mu}{\B(\mu,\nu)}
\dfrac{e^{\mu({\phi-\log(2)})}}{(1+2e^{\phi-\log(2)})^{\mu+\nu}} =
\dfrac{1}{\B(\mu,\nu)}
\dfrac{e^{\mu\phi}}{(1+e^\phi)^{\mu+\nu}}.\qedhere
\end{align*}
\end{proof}

\noindent
{\bf Proposition A.2.} {\it
Let $\gamma$ have density $f(\gamma) =
\dfrac{2^{\mu}}{\B(\mu,\nu)}\dfrac{\gamma^{\mu-1} \;(1-\gamma)^{\nu-1}}{
(1+\gamma)^{\mu+\nu}}, \; 0<\gamma<1$. Then, $\pi=\dfrac{1-\gamma
}{1+\gamma}\sim\Be(\nu,\mu)$.}
\begin{proof}
We have $\gamma= \dfrac{1-\pi}{1+\pi}$. So
$1-\gamma= \dfrac{2\pi}{1+\pi}$, $1+\gamma= \dfrac{2}{1+\pi}$, and
$\dfrac{d\gamma}{d\pi} = \dfrac{-2}{(1+\pi)^2}$. Thus, we have
\[
f(\pi) = \dfrac{2^{\mu}}{\B(\mu,\nu)}
\dfrac{(1-\pi)^{\mu-1}}{(1+\pi)^{\mu-1}}\dfrac{2^{\nu-1}\pi^{\nu
-1}}{(1+\pi)^{\nu-1}}
\dfrac{(1+\pi)^{\mu+\nu}}{2^{\mu+\nu}}\dfrac{2}{(1+\pi)^2}=
\dfrac{1}{\B(\nu, \mu)}\pi^{\nu-1}(1-\pi)^{\mu-1}.\qedhere
\]
\end{proof}

The marginal prior distribution
$\pi(\gamma)$ has the following properties. When $\mu=\nu=1$, $\pi
(\gamma)$ is decreasing. When
$\mu=1$ and $\nu<1$, $\pi(\gamma)$ is U-shaped, the anti-mode being at
$\nu$.
When $\mu=1$ and $\nu>1$ or $\nu=1$ and $\mu<1$ or $\mu>1$ and $\nu<1$,
$\pi(\gamma)$ is decreasing. When $\nu=1$ and $1<\mu<3$, $\pi(\gamma)$
is unimodal
and the mode is at $\gamma=\dfrac{\mu-1}{2}$. When $\nu=1$ and $\mu\geq3$,
it is J-shaped. When $\mu,\nu>1$, $\pi(\gamma)$ is unimodal and the
mode is at
$\gamma= \dfrac{2(\mu-1)} {2\nu+\mu-1+\sqrt{\Lambda}}$,
where $\Lambda=(2\nu+\mu-1)^2-8(\mu-1)$. When $\mu<1$ and $\nu<1$,
$\pi(\gamma)$ is U-shaped, and the anti-mode is at $\gamma = \dfrac
{2\nu+\mu-1+\sqrt{\Lambda}}{4}$.

\section*{Appendix 4: Posterior Distribution of $\Delta$ and $R$}\label{sec4app}

To derive the posterior distribution of $\Delta$ and $R$, we consider
the parameterization $(\gamma, U, V)$.
Under this parameterization, the posterior distributions of interest
belong to the family%\vadjust{\eject}
\begin{eqnarray*}
\pi\big(\gamma, u, v \,|D\big) & = & \dfrac{1}{K} \;\dfrac
{2^{m_{1+}+\alpha}}{\B(m_{1+}+\alpha, m_{2+}+\beta)}
\dfrac{\gamma^{m_{1+}+\alpha-1}\;(1-\gamma)^{m_{2+}+\beta-1}}{(1+\gamma
)^{m_{1+}+m_{2+}+\alpha+\beta}}\;(u+rv)^{d}\\
& & \times\dfrac{u^{m_{10}+m_{20}+a_0-1}(1-u)^{m_{00}+b_0-1}}{\B
(m_{10}+m_{20}+a_0, m_{00}+b_0)}\dfrac{v^{m_{11}+m_{21}+a_1-1}\;
(1-v)^{m_{01}+b_1-1}}{\B(m_{11}+m_{21}+a_1, m_{01}+b_1)}\,,
\end{eqnarray*}
where $0<\gamma, u, v<1$ and $K$ is the normalizing constant. The
choice $d=1/2,a_0=b_0=a_1=b_1=1/2$
corresponds to Jeffreys' posterior distribution and the choice
$d=0,a_0=b_0=a_1=b_1=1/2$ corresponds
to Bernardo's posterior distribution.

\noindent
{\bf Proposition A.3.} {\it
The posterior distribution of the risk difference, $\Delta=\lambda
_1-\lambda_0=\dfrac{V-U}{1+\gamma}$,
has the complex integral form
\begin{align*}
& \; \pi\big(\Delta\,|D\big) \\
= & \; \dfrac{1}{K} \;\int_{-\frac{1+|\Delta|}{|\Delta|}}^{\frac
{1-|\Delta|}
{|\Delta|}} \dfrac{2^{m_{1+}+\alpha}}{\B(m_{1+}+\alpha, m_{2+}+\beta)}
\dfrac{\gamma^{m_{1+}+\alpha-1}
(1-\gamma)^{m_{2+}+\beta-1}}{(1+\gamma)^{m_{1+}+m_{2+}+\alpha+\beta}} \\
& \times\int_{\max(0,-(1+\gamma)\Delta)}^{\min(1,1-(1+\gamma)\Delta
)}\Big((1+r)u+r(1+\gamma)\Delta\Big)^{d}
\dfrac{u^{m_{10}+m_{20}+a_0-1}(1-u)^{m_{00}+b_0-1}}{\B
(m_{10}+m_{20}+a_0, m_{00}+b_0)}\nonumber\\
& \times\dfrac{\Big(u+(1+\gamma)\Delta\Big)^{m_{11}+m_{21}+a_1-1}
\Big(1-u-(1+\gamma)\Delta\Big)^{m_{01}+b_1-1}}{\B(m_{11}+m_{21}+a_1,
m_{01}+b_1)}\,du\,d\gamma.
\end{align*}
We also have
\begin{align*}
& P\big(\lambda_1-\lambda_0>\Delta_0|D\big) \\
= & \;
\dfrac{2^{m_{1+}+\frac{1}{2}}}{\B\big(m_{1+}+\frac{1}{2},m_{2+}+\frac
{1}{2}\big)}
\int_{0}^1 \int_0^1 \int_{v>u+\Delta_0/(1+\gamma)}^1 \dfrac{\gamma
^{m_{1+}+\frac{1}{2}-1}
(1-\gamma)^{m_{2+}+\frac{1}{2}-1}}{(1+\gamma)^{m_{1+}+m_{2+}+1}}\\
& \times\dfrac{(u+rv)^{d}}{K} \dfrac{u^{m_{10}+m_{20}+\frac
{1}{2}-1}(1-u)^{m_{00}+\frac{1}{2}-1}}{\B\big(m_{10}+m_{20}+\frac
{1}{2},m_{00}+\frac{1}{2}\big)} \dfrac{v^{m_{11}+m_{21}+\frac{1}{2}-1}
(1-v)^{m_{01}+\frac{1}{2}-1}}{\B\big(m_{11}+m_{21}+\frac
{1}{2},m_{01}+\frac{1}{2}\big)} dv du d\gamma.
\end{align*}
When $\Delta_0=0$, this expression depends no longer on $\gamma$ and
when $d=0$ it is even simpler.
}

\noindent
{\bf Proposition A.4.} {\it
The posterior distribution of the risk ratio, $R=\dfrac{\lambda
_1}{\lambda_0}=\dfrac{V}{U}$, does not depend on
$\gamma$ and it has the simpler integral form
\begin{align*}
\pi\big(R\,|D\big) = \left\{
\begin{array}{rr}
\displaystyle{\frac{R^{m_{11}+m_{21}+a_1-1}(1+rR)^d}{K\,K^\ast} \int_0^1
u^{m_{1+}+m_{2+}+d+a_0+a_1-1}(1-u)^{m_{00}+b_0-1}} \\
\qquad\displaystyle{\times\left(1-Ru\right)^{m_{01}+b_1-1} \,du}\,
,\qquad0<R \leq1,\\
\displaystyle{\frac{R^{-(m_{10}+m_{20}+a_1+d+1)}(1+rR)^d}{K\,K^\ast}
\int_0^1
v^{m_{1+}+m_{2+}+d+a_0+a_1-1}(1-v)^{m_{01}+b_1-1}}\\
\qquad\displaystyle{\times\left(1-Rv\right)^{m_{00}+b_0-1} \,dv}\,
,\qquad\qquad R>1,
\end{array}
\right.
\end{align*}
where $K^\ast= \B(m_{10}+m_{20}+a_0, m_{00}+b_0)\B(m_{11}+m_{21}+a_1,
m_{01}+b_1)$. When $d=0$,
it can be shown that $\pi\big(R\,|D\big)$ is unimodal when $m_{00}$ and
$m_{01}$ are positive.
We also have
\begin{align*}
P\big(R>R_0|D\big)
= \int_0^1 \int_{v>R_0 u}^1 \Big[ & \dfrac{(u+rv)^{d}}{K}
\dfrac{u^{m_{10}+m_{20}+\frac{1}{2}-1}(1-u)^{m_{00}+\frac{1}{2}-1}}{\B
\big(m_{10}+m_{20}+\frac{1}{2},m_{00}+\frac{1}{2}\big)} \\
& \times
\dfrac{v^{m_{11}+m_{21}+\frac{1}{2}-1}
(1-v)^{m_{01}+\frac{1}{2}-1}}{\B\big(m_{11}+m_{21}+\frac
{1}{2},m_{01}+\frac{1}{2}\big)} \Big] dv du.
\end{align*}
}

%
% References
%
\end{document}